%% file: TD-GCM_revision.tex
\documentclass{article}

\usepackage{arxiv}

\usepackage[utf8]{inputenc} 
\usepackage[T1]{fontenc}    
\usepackage{hyperref}       
\usepackage{url}            
\usepackage{booktabs}       
\usepackage{amsfonts}       
\usepackage{nicefrac}       
\usepackage{microtype}      
\usepackage{lipsum}         
\usepackage{graphicx}
\usepackage{natbib}
\usepackage{doi}

\DeclareUnicodeCharacter{05F3}{`}

\usepackage{times}

\usepackage[english]{babel}
\usepackage{amsmath}
\usepackage{amssymb}

\usepackage{xcolor}
\usepackage{enumitem}
\setenumerate[1]{label=$\bullet$}
\setenumerate[2]{label=$\circ$}

\usepackage{hyperref}
\hypersetup{
	colorlinks=true,
	citecolor=blue,
	linkcolor=blue
}

\usepackage{subcaption}
\usepackage{booktabs}

\usepackage{varwidth}

\usepackage{tikz}
\usetikzlibrary{decorations.pathreplacing,calligraphy}
\usetikzlibrary{shapes.geometric}
\usepackage{pgfplots}
\usetikzlibrary{pgfplots.groupplots}
\pgfplotsset{compat=1.16}
\usetikzlibrary{positioning}
\usetikzlibrary{backgrounds}

\usepackage{marvosym}

\usepackage{diagbox}

\definecolor{limegreen}{RGB}{50,205,50}
\colorlet{mygreen}{green!60!black}

\def\profiles{speed-profiles}
\def\profile{speed-profile}
\def\Profiles{Speed-profiles}

\usepackage{etoolbox}
\def\refEqu#1{(\ref{#1})}
\newcommand{\refFigb}[2][]{\ifblank{#1}{(Fig. #2)}{(Fig. #2 #1)}}
\newcommand{\refFigcf}[2][]{\ifblank{#1}{(cf. Fig. #2)}{(cf. Fig. #2 #1)}}
\newcommand{\refFig}[2][]{\ifblank{#1}{Figure #2}{Figure #2 #1}}
\newcommand{\refTabb}[2][]{\ifblank{#1}{(Tab. #2)}{(Tab. #2 #1)}}
\newcommand{\refTabcf}[2][]{\ifblank{#1}{(cf. Tab. #2)}{(cf. Tab. #2 #1)}}
\newcommand{\refTab}[2][]{\ifblank{#1}{Table #2}{Table #2 #1}}
\def\refSec#1{(Sect. \ref{#1})}

\def\ALGO{PMC-CG}

\def\set#1#2{\left\{\left. {#1} \, \right|\, {#2} \right\}}

\def\zaK{z_a}
\def\zaH{z_a}

\def\halt#1{h(t_{#1},t_{{#1}^\prime})}

\title{Solving the Real-Time Train Dispatching Problem by Column Generation}

\date{January 22, 2024}

\author{ \href{https://orcid.org/0009-0009-9925-3721}{\includegraphics[scale=0.06]{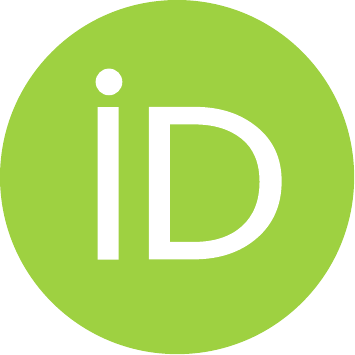}}\hspace{1mm}Maik Schälicke \thanks{maik.schaelicke1@tu-dresden.de}\\
	Dresden University of Technology \\
	Chair of Traffic Flow Sciene\\
	Hettnerstraße 1, 01069 Dresden, Germany \\
	\texttt{https://tu-dresden.de/bu/verkehr/ila/vkstrl}
	\And
	Karl Nachtigall \\
	Dresden University of Technology \\
	Chair of Traffic Flow Sciene\\
	Hettnerstraße 1, 01069 Dresden, Germany 
}

\renewcommand{\headeright}{}
\renewcommand{\undertitle}{}
\renewcommand{\shorttitle}{Solving the TDP by CG}

\hypersetup{
pdftitle={\ALGO},
pdfauthor={Maik Schälicke},
pdfkeywords={Train Dispatching, Column Generation, Speed Profiles},
}

\begin{document}
\maketitle

\begin{abstract}
Disruptions in the operational flow of rail traffic can lead to conflicts between train movements, such that a scheduled timetable can no longer be realised. This is where dispatching is applied, existing conflicts are resolved and a dispatching timetable is provided. In the process, train paths are varied in their spatio-temporal course. This is called the train dispatching problem (TDP), which consists of selecting conflict-free train paths with minimum delay. Starting from a path-oriented formulation of the TDP, a binary linear decision model is introduced. For each possible train path, a binary decision variable indicates whether the train path is used by the request or not. Such a train path is constructed from a set of predefined path parts (\profiles{}) within a time-space network. Instead of modelling pairwise conflicts, stronger MIP formulation are achieved by a cliques formulated over the complete train path. The combinatorics of \profiles{} and different departure times results in a large number of possible train paths, so that the column generation method is used here. Within the subproblem, the shadow prices of conflict cliques must be taken into account. When constructing a new train path, it must be determined whether this train path belongs to a clique or not. This problem is tackled by a MIP. The methodology is tested on instances from a dispatching area in Germany. Numerical results show that the presented method achieves acceptable computation times with good solution quality while meeting the requirements for real-time dispatching.
\end{abstract}


\section{Introduction to the TDP}

In railway operations, train dispatching stands as one of the most crucial and challenging tasks. In the event of operational disruptions, it may become impossible to execute the planned operational schedule. In such instances, dispatchers must contribute to achieving the dispatching objective: the fastest possible return to the planned and conflict-free timetable. Considering all involved train services and the available infrastructure, this poses an exceedingly complex task. To manage the dispatcher workload, they are responsible for a spatially delimited area known as the dispatching area. Moreover, only train services within a temporally dispatching horizon are considered. Within this spatially and temporally delimited area, dispatchers employ dispatching measures such as rerouting or rescheduling to achieve the dispatching objective. By applying these measures, they aim to optimize network capacity utilization to minimize delays for train services. Hence, a train service seeks a conflict-free train path with minimal delays. Thus, the Train Dispatching Problem (TDP) can be defined as the conflict-free selection of train paths for the set of train services within the dispatching area and horizon, aiming to minimize delays.

To further aid dispatchers in solving the TDP, computer-assisted methods are employed. This paper introduces an algorithmic approach to solve the TDP. Section \ref{sec:literature} provides a brief overview of existing approaches. The subsequent Section \ref{sec:solution_algorithm} details the approach presented in this study. Numerical computations and their results are presented in Section \ref{sec:numerical_experiments}. Finally, a concise summary and outlook are provided \refSec{sec:conclusion}.

\section{Related Research and This Work} \label{sec:literature}

According to \cite{Fang2015}, alternative graphs (AG) \citep[][]{Mascis2002} are widely applied, particularly for conflict resolution. In the context of real-time train dispatching, \cite{DAriano2007a} utilise an AG-based formulation. In this approach, the routes are fixed, and the selection of arcs from pairwise incompatible arcs represents a valid scheduling or sequence of trains. To solve their AG model, \citeauthor{DAriano2007a} present a branch-and-bound algorithm. In a subsequent work, the real-time traffic management system ROMA (Railway Traffic Optimization by Means of Alternative Graphs) is introduced \citep{DAriano2008}. The model from \cite{DAriano2007a} is extended to include train routing, using a local search algorithm. For routing and scheduling, \cite{Corman2010a} present a tabu search algorithm. In \cite{Sama2017}, different strategies for exploring neighbourhoods are examined and compared with the local search by \cite{DAriano2008} and the tabu search algorithm by \cite{Corman2010a}. It has been shown that optimal solutions are only available for smaller instances \citep{Corman2010a}, or at least, an advantage in terms of computation time can be achieved \citep{Sama2017}. However, the full potential of train routing cannot be fully realised through AG formulations \citep{DAriano2008}.

To overcome the limitations of AG, mixed-integer programming (MIP) formulations are employed, which find broad applications in real-time train scheduling and dispatching \citep{Fang2015}. In MIP formulations, the time at which a train traverses a track component of its selected route is expressed using a continuous variable. This necessitates disjunctive conditions to decide on the sequence in which trains traverse the track component, and linear constraints are formed using binary variables and big-M constraints. Conflict resolution is achieved explicitly by using minimum headway times based on the blocking time model \citep[][p. 26ff]{Pachl2021}. \cite{Pellegrini2012} present a MIP that includes a detailed modelling of the infrastructure based on track circuits. Test instances within a dispatching area around Lille-Flandres station are optimally solved. In \cite{Pellegrini2014}, the MIP is used within a rolling horizon approach. In a subsequent work by \cite{Pellegrini2015}, a heuristic is introduced for these MIP formulations, allowing for high-quality solutions in real-time train dispatching in complex infrastructure areas. Providing high-quality solutions while meeting real-time requirements is one of the major challenges in real-time applications. Such a support system was applied in Norway, based on a decomposition approach by \cite{Lamorgese2012}. \citeauthor{Lamorgese2012} use a MIP formulation for real-time train dispatching and geographically divide the problem into a line problem and a station problem, reducing the number of big-M constraints. The line problem serves as the master problem, and the station problem forms the subproblem. Based on this decomposition, the solution is obtained similarly to Benders' decomposition, by generating feasible cuts. In \cite{Lamorgese2016}, a heuristic for this decomposition approach was presented, improving the speed of the solution process. Comprehensive overviews of MIP models for real-time train dispatching in railway operations can be found in the review articles by \cite{Cacchiani2014} and \cite{Fang2015}.

The drawback of MIP models lies in their linear programming (LP) relaxation. This relaxation depends on the choice of big-M, but generally represents a weak formulation \citep{Lodi2010, Bonami2015}. In contrast to disjunctive MIP models, time-index models \citep{Dyer1990} discretise the time horizon into time intervals without the need for big-M constraints. The allocation of a track component to a discrete time interval represents a spatial-temporal resource utilisation. For a train path, a path within a space-time network represents an allocation of these resources. In an IP (integer programming) formulation, binary variables are used to decide on the utilisation of these spatial-temporal resources for individual trains. Conflict avoidance is guaranteed through set packing constraints for individual spatial-temporal resources. The critical point in time-index models is the discretisation into time intervals. Finer granularity in the time component results in a large number of binary variables. To manage this, decomposition techniques can be used by splitting the original problem into a master problem and multiple subproblems. This decomposition typically relates to trains. Independent subproblems, representing train-specific issues, can be formulated as a path-based problem \citep{Leutwiler2023}.

In the literature, decomposition approaches for time-index models are primarily found in the formulation of the Train Timetabling Problem (TTP). Similar approaches can be applied to the Train Dispatching Problem (TDP), as it is an online manifestation of the TTP. However, limitations arise due to the macroscopic perspective in the TTP. The basic idea behind decomposition is to formulate the feasibility conditions, i.e., the requirement for conflict avoidance, in a master problem. To achieve this, cliques of pairwise incompatible arcs in a space-time network, representing the utilized track components of a train path, are formed and expressed as set packing constraints in the integer programming (IP) formulation \citep{Brannlund1998, Cacchiani2012}. \cite{Caprara2002}, \cite{Caprara2006}, and \cite{Cacchiani2008} formulate capacity constraints over the nodes of the space-time network instead of the arcs to reduce the number of constraints. Another approach is pursued in the works of \cite{Borndorfer2007} and \cite{Borndorfer2010}, where feasible arc configurations are assigned to each train path, forming a feasible set of arcs, i.e., conflict free track components. In this approach, the focus is not on conflict exclusion but rather on selecting feasible routes. Depending on the formulation of the master problem, in all cases, train paths are generated in the subproblems, which are obtained by applying Lagrangian relaxation \citep{Brannlund1998, Caprara2002, Caprara2006, Cacchiani2012} or through column generation (CG) \citep{Borndorfer2007, Borndorfer2010, Cacchiani2008}.

In the context of real-time train dispatching, \cite{Lusby2013} present a path-based formulation. In an integer programming (IP) model, a binary variable is used to assign a train path to each train service. \citeauthor{Lusby2013} consider three factors for possible paths: crossing opportunities, train kinematics, and the signal system. Considering these factors, all possible paths, i.e., train paths, are formulated in a tree structure. For track capacities, the spatial-temporal resources used by the path are considered. Using set packing formulations, each resource can be claimed by only one train service. The problem is decomposed similarly to column generation. However, unlike classical column generation, the subproblem is not formulated as an optimization problem. Instead, pricing and generating new train paths are done by recursively traversing the developed tree structure. A branch-and-price algorithm is used for generating integer solutions and evaluating their quality.

In \cite{Meng2014}, a decomposition approach is applied to an edge-based IP. The capacity constraint is formulated for train pairs over individual track components. These constraints affect multiple trains and are identified as complicating constraints. Applying Lagrangian relaxation to these capacity constraints results in train-specific subproblems in the form of time-dependent shortest path problems. \citeauthor{Meng2014} present a label-correcting algorithm for solving these subproblems.

In \cite{Bettinelli2017}, the problem is divided into a construction phase and a shaking phase. In the construction phase, a train path is generated by solving a shortest path problem within a space-time network individually for each train service. In the shaking phase, the train order is varied while considering dispatching rules and conflicts. \citeauthor{Bettinelli2017} use and compare a local search algorithm called reduced variable neighbourhood search and a tabu search for this purpose. The construction and shaking phases are iteratively executed to improve the solution.

For real-time train dispatching, \cite{Reynolds2020} and \cite{Reynolds2022} use an edge-based IP. They restrict the simultaneous usage of track components in a way that they can be occupied or banned. The IP formulation of these capacity constraints is based on set packing conditions. Using Danzig-Wolfe decomposition, the authors transition from edge-based to path-based formulation and solve it using a branch-and-price algorithm. Subproblems can also be solved using shortest path algorithms in this approach.

For a comprehensive overview of decomposition approaches for railway scheduling problems, please refer to \cite{Leutwiler2023}.

In order to be able to react at any time online on modified traffic situations, the calculation time to solve the problems should be in the range of seconds. For these reasons, we do not determine the global optimum with a branch-and-price method, but only solve the relaxation with column generation from which an integer solution is determined at the end. By means of a Lagrange relaxation we can specify lower limits to the global optimum, which in many cases show a good quality of the calculated solution. Traditionally, column generation is aimed to solve the LP-relaxation as quickly as possible without any concern for the integer properties of the columns formed. In our approach, we aim to generate columns, which will forming a good integer solution. Instead of using pairwise conflict constraints, we will use maximum conflict cliques, which  is a much stronger model formulation. Solving the pricing problem with maximum cliques as conflict restrictions is a big challenge, because optimizing a new train path must take all shadow prices of the cliques that contain this train path into account.

\section{TDP Solution Algorithm} \label{sec:solution_algorithm}

\begin{figure}[ht!]
	\centering
	\input{img/overview}
	\caption{Two-part approach for the TPD with an independent real-time train dispatching algorithm decoupled from the train path generation.}
	\label{fig:method}
\end{figure}
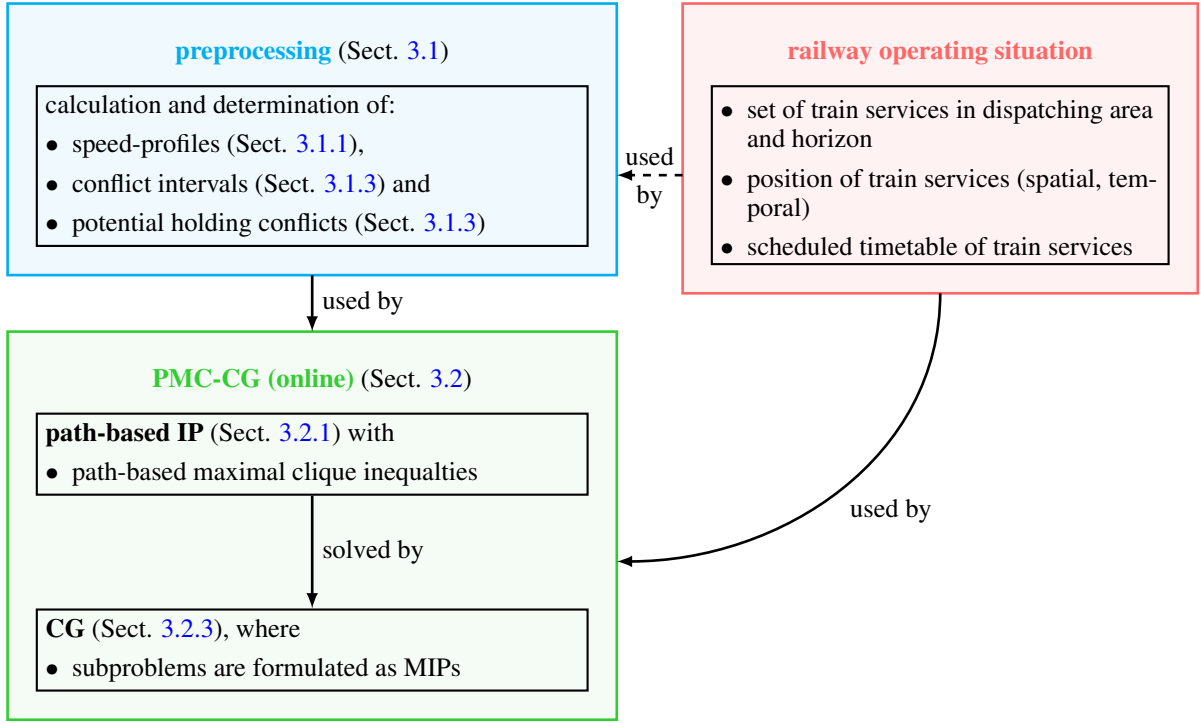

A train path $\left( \left( v_1, t_{v_1} \right), ..., \left( v_n, t_{v_n}  \right) \right)$ will be described by sequence of path segments $\left( v_i, t_{v_i} \right)$, so called \profiles{} $v_i$ and time points $t_{v_i}$. A \profile{} $v$ has only a temporal length, but no fixed time. Minimum headway times between train paths will be calculated from time points $t_v$ and the minimum headway times between \profiles{}. The driving dynamic calculation of the \profiles{} and their minimum headway times are carried out as a preprocess and thus, strictly separated from the optimization core \refFigb{\ref{fig:method}}.

This approach has the advantage, that the optimization is independent of the \profile{} construction process. \Profiles{} can be generated using commercial software such as LUKS \citep{Luks} or through custom-developed methods based on space-time networks \citep[e.g.][]{Lusby2013}. The preprocessing part can occur in a more strategic planning manner, that means a set of \profiles{} can be generated beforehand and stored in a database. The database needs to be updated to reflect the current railway traffic situation.

The optimization core is a path-based formulation, where each train path for every train service is modelled by a binary decision variable. Conflict restrictions between train paths are formulated by maximum cliques. Due to the huge number of variables we will use a column generation algorithm (\ALGO{}) to solve the TDP.


\subsection{Preprocessing} \label{sec:preprocessing}

\subsubsection{Speed-Profiles} \label{sec:speed_profiles}

The concept of \profiles{} corresponds to snippets described in \cite{Nachtigall2016} or \cite{Dahms2019}, used for generating schedules for freight trains by providing train paths for these trains. Mathematically, a \profile{} represents a space-time function. The spatial extent for a \profile{} is limited by points where dispatching decisions can be made, referred to as dispatching points. The \profiles{} given by the preprocess define the potential of possible dispatching decisions. Regardless of the optimization, they can be designed from simple to arbitrarily complex.

For the sequence of block sections utilized by a \profile{} $v$, the calculation of running time is conducted based on a microscopically accurate vehicle dynamics model.
\refFig{\ref{fig:speed_profile_and_train_path}} contains an example of constructing train paths for a train service $r$, depicting the scheduled route from Station $S_1$ via $S_2$ to $S_3$. In \refFig[(a)]{\ref{fig:speed_profile_and_train_path}}, the utilized block sections are illustrated for the color-highlighted \profiles{}. For example, the \profile{} $v_1$ possesses $\left( b_{10}, b_{11}, b_{13}, b_{15}, b_{16}, b_{18}, b_{21} \right)$ as the utilized sequence of block sections. As marginal conditions, it is specified that the initial and final velocities are set to zero. Adhering to these marginal conditions and considering the principles of vehicle dynamics, acceleration occurs to reach a predetermined maximum speed, followed by maintaining a constant speed and deceleration to return to zero velocity. When defining the maximum speed, it is crucial to note the feasibility of decelerating to reach zero velocity. This results in a running time, exemplified by a \profile{} $v_1$ as shown in \refFig[(b)]{\ref{fig:speed_profile_and_train_path}}. This running time also includes the scheduled dwell time at the dispatching point where a \profile{} arrives.

For variations concerning the \profiles{}, the following possibilities are considered:
\begin{enumerate}[label=(V\arabic*)]
	\item Varying the predetermined maximum speed while maintaining an unaltered sequence of block sections.  \label{variation1}
	\item Altering the sequence of block sections while maintaining a predetermined maximum speed. \label{variation2}
	\item A combination of varying both the sequence of block sections and the maximum speed. \label{variation3}
\end{enumerate}
For speed variation \ref{variation1}, the initial choice might be a nominal speed predetermined by the timetable. To compensate for delays, a \profile{} with a shorter travel time is generated. Conversely, to avoid catching up to a preceding train, a \profile{} with a longer travel time is generated. In \refFig[(a)]{\ref{fig:speed_profile_and_train_path}}, the \profiles{} $v_2$ and $v_3$ utilize the same infrastructure, yet due to different maximum speeds, they exhibit varying space-time curves and consequently different travel times \refFigb[(b)]{\ref{fig:speed_profile_and_train_path}}. Regarding route variation \ref{variation2}, \profiles{} are initially created for the scheduled train path of the train service, ensuring that following the schedule is always possible. For this purpose, the scheduled platforms for halts are considered as dispatching points, between which \profiles{} are then generated. The route variation process is performed microscopically by deviating from the schedule to approach a different platform during a scheduled halt. An example of this is illustrated by the \profiles{} $v_1$ and $v_4$ in Figure \ref{fig:speed_profile_and_train_path} (a). Additionally, a macroscopic route variation involves deviating from the schedule between two consecutive, scheduled stations $S_1$ and $S_2$ by making one or several additional halts. For an additional halt $X$, it results in a macroscopic travel sequence $\left( S_1, X, S_2 \right)$, where a microscopic route variation occurs again between $\left( S_1, X \right)$ and $\left( X, S_2 \right)$. Macroscopic route variation between two consecutive, scheduled stations $S_1$ and $S_2$ only occurs if the distance covered through the additional halts is not greater than the distance between $S_1$ and $S_2$ multiplied by a detour factor $\rho = 2.5$. This method generates a set of \profiles{} $\mathcal{V} \left( r \right)$ for a train service $r$.

\subsubsection{Train Path Construction from Speed-Profiles} \label{sec:train_path_formulation}

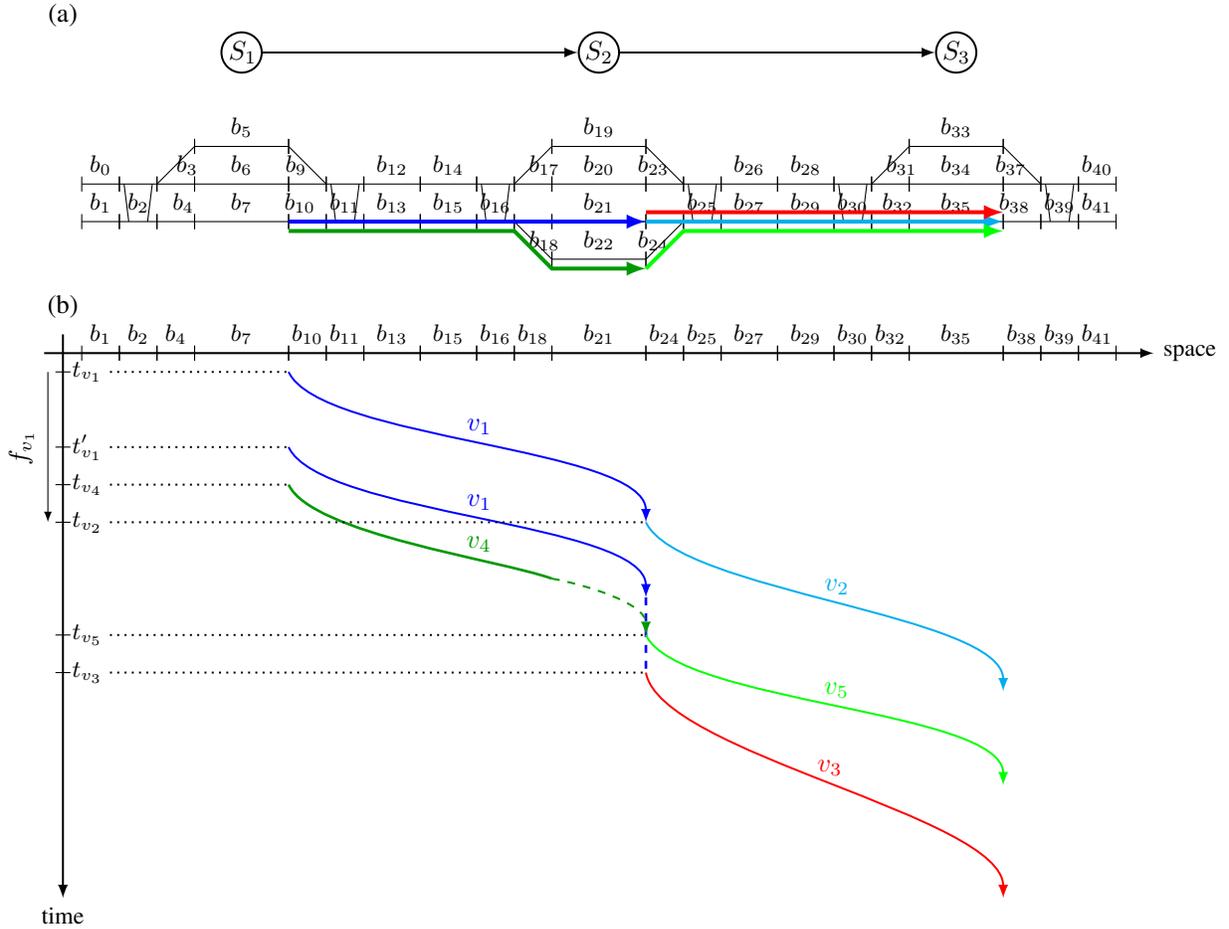
\begin{figure}[ht!]
	\centering
	\input{img/speed_profile_and_train_path}
	\caption{Possible train paths for a train service $r$ from $S_1$ via $S_2$ to $S_3$. (a) Block sections utilized by the color-highlighted \profiles{} $v_1, ..., v_5$. (b) Space-time curve of the \profiles{} and connection of these to train paths $a_1 = \left( \left( v_1, t_{v_1} \right), \left( v_2, t_{v_2} \right) \right)$, $a_2 = \left( \left( v_1, t_{v_1} \right), \left( v_3, t_{v_3} \right) \right)$ and $a_3 = \left( \left( v_4, t_{v_4} \right), \left( v_5, t_{v_5} \right) \right)$.
}
	\label{fig:speed_profile_and_train_path}
\end{figure}

Two \profiles{} $v, v^\prime$ can be connected if, from an infrastructural perspective, there is a continuously passable sequence of block sections. Then, for a profile $v$ with $\mathcal{V}(v)$, the set of possible successor \profiles{} is defined. Each $v$ is assigned a departure time $t_v$, resulting in a pair $(v, t_v)$. A train path $a$ is then a sequence $\left((v_i, t_{v_i})\right)_{i \in \mathbb{N}}$ of these pairs, while adhering to the successor relationships.
In \refFig[(b)]{\ref{fig:speed_profile_and_train_path}}, three possible train paths are depicted for the illustrated train service. The train paths $a_1 = \left( \left( v_1, t_{v_1} \right), \left( v_2, t_{v_2} \right) \right)$ and $a_2 = \left( \left( v_1, t_{v_1}^\prime \right), \left( v_3, t_{v_3} \right) \right)$ differ firstly in the choice of the successor \profile{} of $v_1$. Additionally, $a_2$ utilizes a later time point for the departure time of $v_1$. Deviating from the other depicted train paths in \refFig[(b)]{\ref{fig:speed_profile_and_train_path}}, an additional halt at $S_2$ is planned for $a_2$. With a microscopic rerouting, the train path $a_3 = \left( \left( v_4, t_{v_4} \right), \left( v_5, t_{v_5} \right) \right)$ arrives at a different platform.

A train path for a train service $r$ always starts with a start \profile{} and ends with an end \profile{}. The sets for these distinguished \profiles{} are denoted as $\mathcal{V}_{\text{start}}(r)$ and $\mathcal{V}_{\text{end}}(r)$. During the dispatching, start \profiles{} are those \profiles{} that can be chosen when entering the dispatching area. Similarly, end \profiles{} are those chosen when leaving the dispatching area. 

\subsubsection{Conflict Detection} \label{sec:conflict_detection}

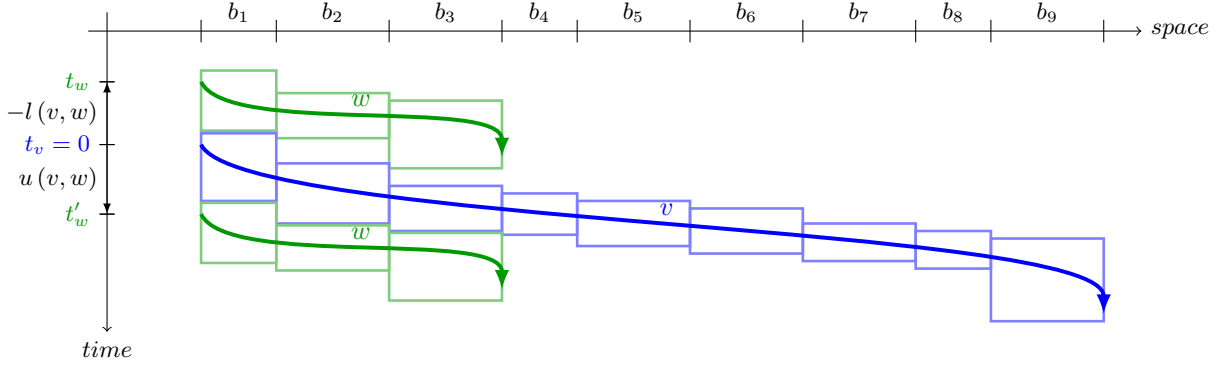
\begin{figure}
	\centering
	\input{img/conflict_1}
	\caption{Conflict interval for a \profile{} $w$ defined by a \profile{} $v$ through the minimum headway times $l \left( v, w \right)$ and $u \left( v, w \right)$.}
	\label{fig:conflict_1}
\end{figure}

For two train paths $a_1$ and $a_2$ with their parts $(v, t_v)$ and $(w, t_w)$ respectively, there is a potential conflict if $v$ and $w$ share a common block section. $(v, t_v)$ and $(w, t_w)$ are conflict-free if and only if 
\begin{equation}
	t_w > t_v + u(v, w) \text{ or } t_w < t_v - l(v, w) \label{equ:disjoint_conflict_condition}
\end{equation}
holds. Here, $u(v, w)$ is the minimum headway time between $v$ and $w$ when $w$ follows $v$, and $l(v, w)$ is the minimum headway time between $v$ and $w$ when $w$ precedes $v$. With $-l(v, w)$ and $u(v, w)$, a conflict interval $K(v, w) = [-l(v, w), u(v, w)]$ is defined, and it follows:
\begin{equation}
	t_w - t_v \in K(v, w) \iff \text{conflict between } v \text{ and } w \label{equ:conflict_condition}
\end{equation}
In \refFig{\ref{fig:conflict_1}}, a conflict between two \profiles{} $v$ and $w$ is avoided by choosing a departure time $t_w$ for $w$ such that $w$ departs before $v$. However, with the choice of $t^\prime_w$, $w$ could also depart conflict-free after $v$. Obviously, in this example, $t_w, t^\prime_w \notin K(v, w)$.

In order to avoid conflicts of overtaking or crossing train paths with halting train paths using the same infrastructure, the halting time 
must be taken into account by enlarging the upper bound of the conflict interval for the planned halting time. 
This halting time $\halt{w} = t_{w^\prime}-t_w-f_w$  depends on the succeeding speed-profile $w^\prime$ and its point of time $t_{w^\prime}.$ 
A train path $a_1$ containing $(v,t_v)$ and a train path $a_2$ containing the part $(w,t_w),(w^\prime,t_{w^\prime})$ 
will have a conflict, if and only if 
\begin{equation}
	t_v - t_w \in [-l(w,v),u(w,v)+\halt{w}] \label{equ:disjoint_halting_conflict_condition_1}
\end{equation}
or
\begin{equation}
	t_v - t_{w^\prime} \in K(w^\prime,v). \label{equ:disjoint_halting_conflict_condition_2}
\end{equation}
For \refEqu{equ:disjoint_halting_conflict_condition_1}, it follows:
\begin{eqnarray}
	& & t_v - t_w \in \left[ -l(w, v), u(w, v) + \halt{w} \right] \nonumber \\
	& \iff & t_v - t_w \ge -l(w, v) \text{ and } t_v - t_{w^\prime} \le u(w, v) - f_w \label{equ:disjoint_halting_conflict_condition_3}
\end{eqnarray}
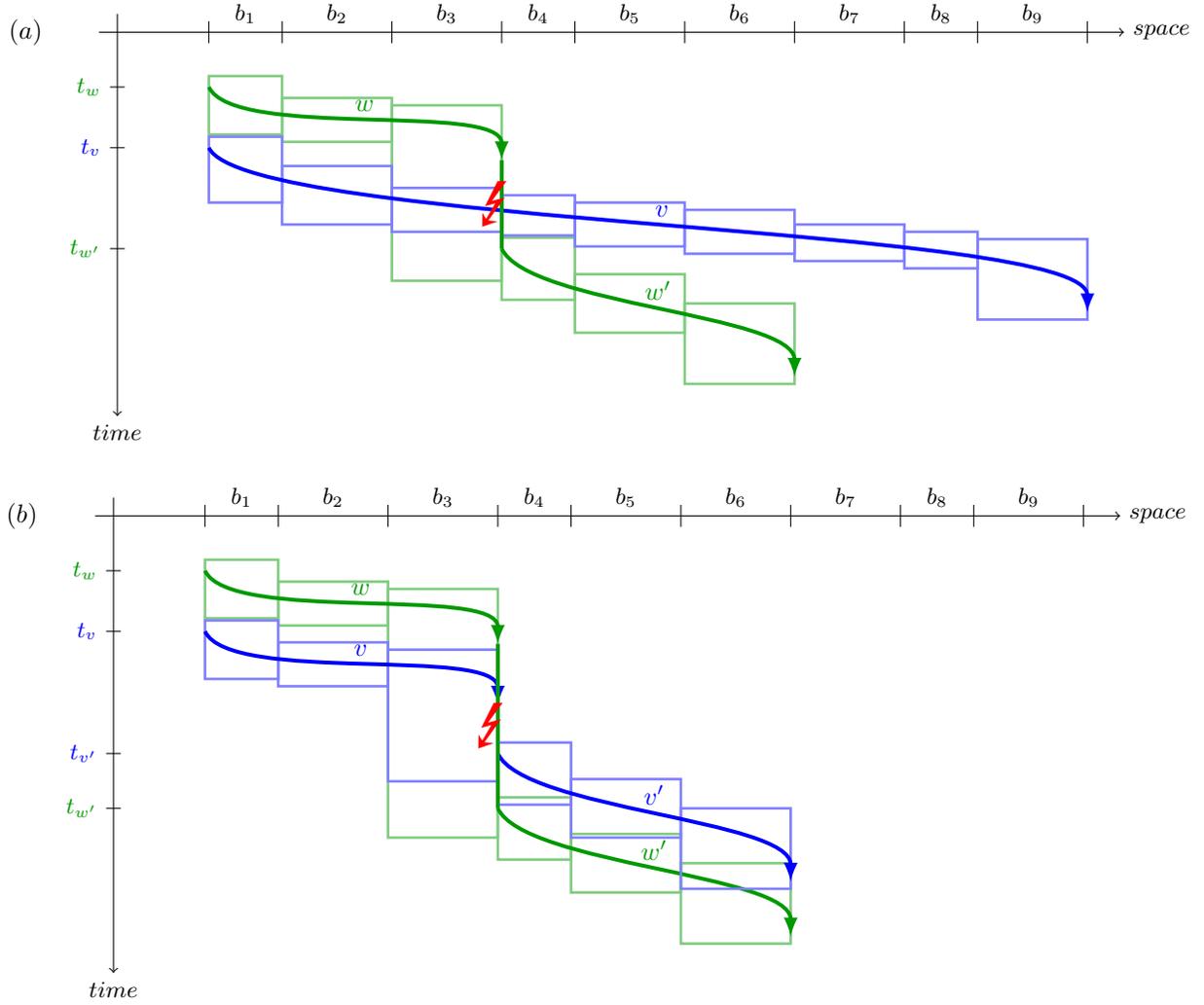
\begin{figure}
	\centering	
	\input{img/conflict_3}
	\input{img/conflict_5}
	\caption{(a) A halting conflict at platform $b_3$ between a crossing and halting train path with parts $\left( v,t_v \right)$ and $\left( w,t_w \right), \left( w^\prime,t_{w^\prime} \right)$, respectively. (b) A halting conflict at platform $b_3$ between two halting train paths with parts $\left( v,t_v \right), \left( v^\prime,t_{v^\prime} \right)$ and $\left( w,t_w \right), \left( w^\prime,t_{w^\prime} \right)$, respectively.}
	\label{fig:conflict_2}
\end{figure}
In \refFig[(a)]{\ref{fig:conflict_2}}, an example of a halting conflict with a crossing \profile{} $v$ is illustrated. In this case, $v$ traverses block $b_3$, which is occupied by the halting of the train path shown in green. Thus, condition \refEqu{equ:disjoint_halting_conflict_condition_1} is satisfied. Condition \refEqu{equ:disjoint_halting_conflict_condition_2} is not fulfilled since there is no conflict with $w^\prime$. Note, that a halting conflict condition $H(v,w,w^\prime)$  is described by one speed-profile $v$ of train path $a_1$ and two suceeding speed-profiles $w,w^\prime$ of train service $a_2$ and its time points $t_v,t_w,t_{w^\prime}$.

A halting conflict can also occur between two halting train paths. Let $a_1$ and $a_2$ be train paths with parts $\left( v, t_v \right), \left( v^\prime, t_{v^\prime} \right)$ and $\left( w, t_w \right), \left( w^\prime, t_{w^\prime} \right)$, respectively. In this case, it is sufficient to consider only the halting time $\halt{w}$ . However, the departure time $t_{v^\prime}$ must be taken into account, and condition \refEqu{equ:disjoint_halting_conflict_condition_2} has to be modified. Then, there is a halting conflict between two halting trains if and only if \refEqu{equ:disjoint_halting_conflict_condition_1} and
\begin{equation}
	t_{v^\prime} - t_{w^\prime} \in K(w^\prime,v^\prime) \label{equ:disjoint_halting_conflict_condition_4}
\end{equation}
holds. The example in \refFig[(b)]{\ref{fig:conflict_2}} illustrates a conflict situation in the case of two halting train paths. In this case, conditions \refEqu{equ:disjoint_halting_conflict_condition_1} and \refEqu{equ:disjoint_halting_conflict_condition_4} are satisfied, resulting in the detection of a conflict. The halting conflict condition $H(v,v^\prime,w,w^\prime)$ for two halting train paths is described by two suceeding \profiles{} $v, v^\prime$ of train path $a_1$ and two suceeding \profiles{} $w,w^\prime$ of train service $a_2$ and its time points $t_v,t_{v^\prime},t_w,t_{w^\prime}$. 

\subsection{\ALGO} \label{sec:PMC-CG}

\subsubsection{Path-oriented IP Formulation} \label{sec:path_based_model}

In the path-based formulation, a path represents a sequence of \profiles{}, in other words, a train path. For each train path, a binary variable is introduced as follows:
\begin{equation*}
	x_a = 
	\begin{cases}
		1 & \text{, if train path $a$ is selceted} \\
		0 & \text{, otherwise} 
	\end{cases}
\end{equation*}
Moreover, $c_a$ denotes the cost associated with train path $a$. In this context, the costs represent the delay incurred by train path $a$ when leaving the dispatching area. The TDP now involves selecting a train path from the set of all possible train paths for each of the train services. Here, the train services considered are those located within the dispatching area and time horizon. This selection of train paths must be conflict-free. 
In this paper we will use column generation for solving the relaxation. When using column generation for a mixed integer problem, the pricing steps are usually performed within the relaxation. Hence, it is very important to have a very strong model formulation. Instead of using pairwise conflict constraints, it is much more stronger to use maximum conflict cliques. A clique naturally arises, if a subset of train paths uses one common block of infrastructure. Those block based cliques are not necessarily maximum. Here, a cross-block conflict consideration takes place. Block-crossing means that cliques are formed by considering conflicts not only within a single block section but across block sections related to the entire train path. Thus, this approach examines network-wide interconnections of conflicts. Train path-based cliques correspond in form to the path-based cliques described in \citet{Caprara2010}, which, as described there, generally represent a stronger LP formulation. The modelling of conflicts is thus carried out using train path-based maximal cliques $C$ within an undirected  conflict graph $G = \left( \mathcal{A}, \mathcal{E} \right)$. In this graph $G$, the set of nodes is represented by the train paths $a \in \mathcal{A}$, and there is an edge $\left( a_1, a_2 \right) \in \mathcal{E}$ if there is a conflict between $a_1$ and $a_2$. The entire conflict situation is described by $\mathcal{C}$, the set of all maximal cliques $C$. In the IP, these cliques are formulated in the form of set-packing constraints. With the objective of the TDP, the following complete model is then derived:
\begin{eqnarray}
	& \displaystyle \sum_{a:a \in \mathcal{A}} c_a x_a \rightarrow \min & \label{equ:ip1} \\
	\forall r \in \mathcal{R}: & \displaystyle \sum_{a:a \in \mathcal{A}\left( r \right)} x_a & = 1 \label{equ:ip2} \\
	\forall C \in \mathcal{C}: & \displaystyle \sum_{a:a \in C} x_a & \leq 1 \label{equ:ip3} \\
	& x_a \in \left\lbrace 0,1 \right\rbrace, a \in \mathcal{A} & \label{equ:ip4}
\end{eqnarray}
The objective function \refEqu{equ:ip1} minimizes the costs of the train paths. Constraint \refEqu{equ:ip2} ensures that exactly one train path is selected for each train service $r \in \mathcal{R}$. This condition is known as the fulfillment condition. With \refEqu{equ:ip3}, it is ensured, that at most one train path is selected from every $C$, ensuring conflict-free selections. Constraint \refEqu{equ:ip4} represents the integrality condition for the variables $x_a$. The set $\mathcal{A}$ contains all possible train paths for all train services $r \in \mathcal{R}$. The number of train paths in $\mathcal{A}$ is very high because for every possible sequence of \profiles{}, there can be a different choice of departure times. 

\subsubsection{Flow of the \ALGO} \label{sec:overview}
The \ALGO{} uses column generation to solve \refEqu{equ:ip1} - \refEqu{equ:ip4}. Column generation searches for new columns for the LP relaxation of the restricted master problem (rRMP), which means to minimise 
\begin{equation}
	z(rRMP) = \min \set{\sum c_a x_a}{(\ref{equ:ip1}) - (\ref{equ:ip3}), x_a \geq 0, a \in \mathcal{A}^\prime \left( \subseteq \mathcal{A} \right)}.
\end{equation}
Note, that the total optimum of TDP is 
\begin{equation}
	z(TDP) = \min \set{\sum c_a x_a}{(\ref{equ:ip1}) - (\ref{equ:ip3}), x_a \in \lbrace 0, 1 \rbrace, a \in \mathcal{A}}.
\end{equation}
If there are train paths $a_{new}$ with negative reduced costs $\rho_{a_{new}}$, the variable $x_{a_{new}}$ associated with $a_{new}$ is added to the rRMP as a new column, expanding the set $\mathcal{A}^\prime$ with $a_{new}$. 
The determination of $\rho_{a_{new}}$ is based on the dual variables $\alpha_r$ and $\beta_C$ related to the fulfillment conditions \refEqu{equ:ip2} and clique inequalities \refEqu{equ:ip3}, respectively.
The reduced costs are computed as follows:
\begin{equation}
	\forall r \in \mathcal{R} : \rho_{a_{new}} = c_{a_{new}} + \sum_{C \in \mathcal{C}:a_{new} \in C} \beta_C - \alpha_r \label{equ:reduced_cost}
\end{equation}
Here, $c_{a_{new}}$ represents the costs of the new train paths, indicating their delay upon leaving the dispatching area. Since $\alpha_r$ remains constant, the subproblem for each train service $r \in \mathcal{R}$ becomes:
\begin{equation}
	\rho^\star (r) = \min \left\lbrace c_{a_{new}} + \sum_{C \in \mathcal{C}:a_{new} \in C} \beta_C \text{ $\left\vert \right. \; a_{new}$ is a train path for train service $r$ } \right\rbrace \label{equ:generic_pricing}
\end{equation}

\begin{figure}[ht!]
	\centering
	\input{img/algorithm_overview}
	\caption{Overview of the \ALGO{}}
	\label{fig:algorithm_overview}
\end{figure}
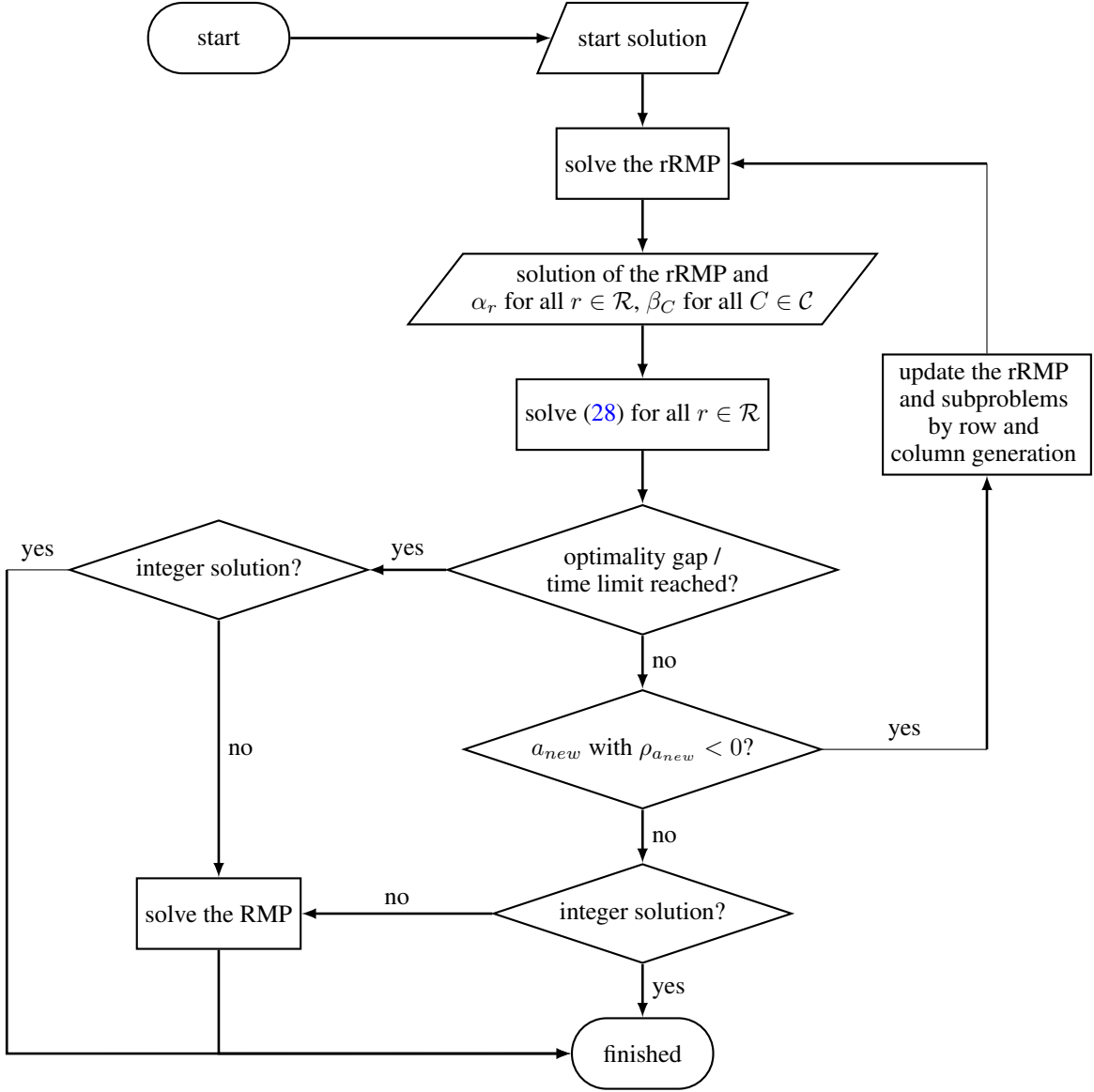

The flow of the PMC-CG is depicted in \refFig{\ref{fig:algorithm_overview}}. The feasible initial solution with the subset $A^\prime\subseteq A$
for column generation is determined using a first-come-first-serve (FCFS) start heuristic. Here, the train
services are sorted based on their arrival time in the dispatching area. Successively, for each train service,
a conflict-free train path is generated added to rRMP\footnote{This is done by solving the subproblem. When a new train path is generated for a train service, new conflict cliques $C$ can arise. If there are conflict cliques, they will be included in the subproblem and the corresponding $\beta_C$ will be set to $\infty$. Then, the subproblem will be solved again. This is carried out, as long as new conflict cliques are detected. If not, there is a conflict free train path for the current train service and it continues with the next train service. This ensures a consistent avoidance of conflicts described by $C$, whereby rerouting and rescheduling also take place.}.

In order to keep the algorithm online, the method terminates, if either a good solution quality or a time limit has been reached.
The quality of the solution is evaluated by the optimality gap, which can be calculated from a lower bound. Following
\cite{Luebbecke2005} and \cite{Luebbecke2011}, Lagrangian relaxation of the conflict constraints \refEqu{equ:ip3} provides
a good lower bound:
\begin{equation}
	lb(rRMP) = z(rRMP) - \sum_{r\in R} \rho^\star(r) \leq z(TDP),
\end{equation}
where $z(rRMP)$ is the actual objective value of the rRMP.
Since our method solves the subproblems by optimality, we are able to compute those values. 

While adding new train path variables to the rRMP the system of clique constraints must be updated.
Some  conflict cliques have to be enlarged and new maximum conflict cliques (= new rows) will arise. This is an essential part of the \ALGO{} and is described in Section \ref{sec:clique_update}.
After finishing the column generation process, a possibly non-integer solution causes to calculate 
the final integer solution of the restricted master problem (RMP).
	
\subsubsection{MIP Subproblem} \label{sec:column_generation}

Upon examining \refEqu{equ:generic_pricing}, it is evident that solving this minimisation problem involves deciding whether $a_{new}$ belongs to a maximal clique $C$ or not. This occurs when $a_{new}$ conflicts with all $a \in C$. The pricing with $\beta_C$ can be avoided by avoiding the conflict with at least one $a \in C$. This can be achieved through rerouting, that is selecting a different sequence of profiles. Alternatively, it could involve rescheduling by adjusting the departure times of \profiles{} or simultaneous rerouting and rescheduling. Typically, these actions result in higher costs $c_{a_{new}}$ due to longer travel times on detours or additional unplanned halts and dwell times. The network-wide conflict consideration in the individual $C \in \mathcal{C}$ now makes it challenging to determine where and how rerouting and rescheduling should occur to achieve the best possible trade-off between costs $c_{a_{new}}$ and $\beta_C$, leading to minimal reduced costs.

This decision is made through a MIP. Let $t_e$ denote the time of the train path $a_{new}$ after using an end \profile{}. For all $v \in \mathcal{V}_{end} \left( r \right)$, $t_e = t_v + f_v$ holds. The end \profile{} ends at a dispatching point, defining the location where the dispatching area is left. Then, $t_e$ is the point of time when the dispatching area is left. The cost $c_{a_{new}}$ results from the difference between $t_e$ and the scheduled time $t_e^{scheduled}$ at this designated dispatching point according to the timetable. A feasible train path is obtained by considering the successor conditions of \profiles{}, which can be realized through flow conservation (Eq. \refEqu{equ:flow1} and \refEqu{equ:flow2}), precedence constraints (Eq. \refEqu{equ:pre1} - \refEqu{equ:pre3}), and the binary variables:
\begin{equation*}
	y_v = 
	\begin{cases}
		1 & \text{, if \profile{} $v$ is selceted} \\
		0 & \text{, otherwise} 	
	\end{cases}
\end{equation*}
The following constraints arise as a result:
\begin{eqnarray}
	& \displaystyle \sum_{v:v \in \mathcal{V}_{start} \left( r \right)} y_v & = 1 \label{equ:flow1} \\
	\forall v \in \mathcal{V} \left( r \right) : & \displaystyle  y_v - \sum_{v^\prime:v^\prime \in \mathcal{V}_{suc} \left( v \right)} y_w & = 0 \label{equ:flow2} \\
	\forall v \in \mathcal{V}_{start} \left( r \right) : & t_v & \geq t_b \label{equ:pre1} \\
	\forall v \in \mathcal{V} \left( r \right), v^\prime \in \mathcal{V}_{suc} \left( v \right) : & t_{v^\prime} & \geq t_v + f_v y_v \label{equ:pre2} \\
	\forall v \in \mathcal{V}_{end} \left( r \right) : & t_e & = t_v + f_vy_v \label{equ:pre3}
\end{eqnarray}
A train path can only start with a starting \profile{}, as ensured by \refEqu{equ:flow1}. With \refEqu{equ:flow2}, exactly one successor must be selected for a \profile{} $v$. Specifically for the starting \profiles{}, \refEqu{equ:pre1} ensures that their departure cannot be before $t_b$, that is the entry time into the dispatching area. If a \profile{} $v$ is chosen, then according to \refEqu{equ:pre2}, for a successor \profile{} $w$, the departure time of $w$ must be greater equal than the departure time of $v$ plus its running time $f_v$. It should be noted that in the case where $t_{v^\prime} > t_v + f_vy_v$ with $y_v = 1$, there is an extended halt beyond the minimum dwell time before the departure of $v^\prime$. Condition \refEqu{equ:pre3} provides the point of time when leaving the dispatching area.

The decision whether $a_{new}$ belongs to a maximal clique $C$ initially involves checking for a conflict with a train path. This process is carried out for all $a \in \bigcup_{C \in \mathcal{C}} C \subseteq \mathcal{A}^\prime$. Following condition \refEqu{equ:conflict_condition} the conflict detection can be made with means of the conflict intervals and binary indicator variables $z_K^l, z_K^u$. This results in the following conditions.
\begin{eqnarray}
	\forall K \left( w, v \right) \in \mathcal{K} : & t_v & \leq t_w - l \left( w , v \right) + M \left( z_K^l - y_v + 1 \right) - \varepsilon \label{equ:drive_up_conflict1} \\
	& t_v & \geq t_w + u \left( w , v \right) - M \left( z_K^u - y_v + 1 \right) + \varepsilon \label{equ:drive_up_conflict2}
\end{eqnarray}
Here, $\mathcal{K}$ is the set of all conflict intervals $K \left( v^\prime, v \right)$ that must be considered for the \profiles{} $v \in \mathcal{V} \left( r \right)$ of the train service $r$. The \profiles{} $w$ are those that are already used by generated train paths $a \in \mathcal{A}^\prime$. The $\varepsilon > 0$ values are sufficiently small so that \refEqu{equ:drive_up_conflict1} and \refEqu{equ:drive_up_conflict2} indeed represent strict less-than and strict greater-than conditions, and $M > 0$ values are sufficiently large. Constraints \refEqu{equ:drive_up_conflict1} and \refEqu{equ:drive_up_conflict2} are active only when the respective \profile{} $v$ is chosen, i.e., $y_v = 1$. With $y_v = 1$, in case of a conflict, i.e., \refEqu{equ:conflict_condition} holds, this implies $z_K^l = z_K^u = 1$. With a binary variable $\zaK$ and
\begin{equation}
\forall a \in \mathcal{A}^\prime : M \zaK \geq \displaystyle \sum_{K: K \in \mathcal{K} \left( a \right)} \left( z_K^l + z_K^u - 1 \right) \label{equ:drive_up_conflict3}
\end{equation}
a conflict with a train path $a$ can be detected. $\mathcal{K} \left( a \right) \subseteq \mathcal{K}$ is the set of conflict intervals that must be considered for the conflict detection between the train path $a$ with the new train path to be generated for $r$. If at least one conflict interval $K \left( v^\prime , v \right) \in \mathcal{K} \left( a \right)$ exists such that $t_v \in K \left( v^\prime , v \right)$, then $\zaK = 1$, indicating the conflict with $a$.

If conditions \refEqu{equ:disjoint_halting_conflict_condition_1}, \refEqu{equ:disjoint_halting_conflict_condition_2} and \refEqu{equ:disjoint_halting_conflict_condition_4} hold, a halting conflict arises. Since condition \refEqu{equ:disjoint_halting_conflict_condition_2} and \refEqu{equ:disjoint_halting_conflict_condition_4} are already included in \refEqu{equ:drive_up_conflict1} - \refEqu{equ:drive_up_conflict3} we must only consider \refEqu{equ:disjoint_halting_conflict_condition_1}. With \refEqu{equ:disjoint_halting_conflict_condition_3} it follows:
\begin{eqnarray}
	\forall H\left( v , w, w^\prime \right) \in \mathcal{H} : & t_v & \leq t_w - l \left( w, v \right) - M \left( z_H^l - y_v + 1 \right) + \varepsilon \label{equ:halting_conflict1} \\
	& t_v & \geq t_{w^\prime} +  u \left( w, v \right) - f_w + M \left( z_H^u - y_v + 1 \right) - \varepsilon \label{equ:halting_conflict2} \\
	\forall a \in \mathcal{A}^\prime : & M z_a & \geq \displaystyle \sum_{H: H \in \mathcal{H} \left( a \right)} \left( z_H^l + z_H^u - 1 \right) \label{equ:halting_conflict3}
\end{eqnarray}
Here, we also use binary indicator variables $z_H^l, z_H^u$ to decide if $t_v - t_w \in \left[ -l(w,v), u(w,v) + h(t_w, t_{w^\prime}) \right]$. This is carried out for every $v \in \mathcal{V} (r)$, where all corresponding $H(v,w,w^\prime)$ will be considered. $\mathcal{H}$ is set of all halting conflict conditions\footnote{It's easy to see that a similar condition holds for the case of two halting trains $H(v,v^\prime,w,w^\prime)$ \refSec{sec:conflict_detection}.}. In \refEqu{equ:halting_conflict1} and \refEqu{equ:halting_conflict2}, appropriately small $\varepsilon$ and large $M$ values must be chosen. If the \profile{} $v$ is not chosen, the constraints \refEqu{equ:halting_conflict1} and \refEqu{equ:halting_conflict2} become inactive. Otherwise a halting conflict is detected if $z_H^l = z_H^u = 1$ and with \refEqu{equ:halting_conflict3} then $\zaH = 1$ and thus a conflict with $a$ follows. For this purpose, constraint \refEqu{equ:halting_conflict3} uses the set $\mathcal{H} \left( a \right) \subseteq \mathcal{H}$, which is the set of halting conflict conditions that must be considered for the potential halting conflicts with $a$. With variables $z_a$ and the constraint
\begin{equation}
	\forall C \in \mathcal{C}: z_C \geq \left( \sum_{a:a \in C \in \mathcal{C}} z_a \right) - \left\lvert C \right\rvert + 1, \label{equ:clique}
\end{equation}
the binary variable $z_C$ determines whether $a_{new}$ belongs to a maximal clique $C$. Here, $\left\lvert C \right\rvert$ represents the cardinality of $C$. The right-hand side of \refEqu{equ:clique} becomes one if conflicts exist with all $a \in C$, consequently resulting in $z_C = 1$. Using $z_C \beta_C$ for all $C \in \mathcal{C}$, the MIP subproblem models the costs for conflicts in the objective function. The resulting complete MIP subproblem for each $r \in \mathcal{R}$ is formulated as follows:
\begin{equation}
	\begin{array}{rc}
		& \displaystyle \left( t_e - t_e^{scheduled} \right) + \sum_{C:C \in \mathcal{C}} z_C \beta_C \rightarrow \min \\ \\
		s.t. & \refEqu{equ:flow1} - \refEqu{equ:clique} \\ \\
		& t_v,t_w,t_e \geq 0, y_v,y_w,z_K^l,z_K^u,z_H^l,z_H^u,z_a,z_C \in \lbrace 0, 1 \rbrace
	\end{array} \label{equ:mip_subproblem}
\end{equation}

\subsubsection{Determining and Updating the Maximal Cliques} \label{sec:clique_update}

An essential component of the \ALGO{} presented here is the determination and updating of the maximal cliques. When new train paths are generated by solving the subproblems \refEqu{equ:mip_subproblem} for all $r \in \mathcal{R}$, the clique inequalities \refEqu{equ:ip3} must be supplemented with the respective variables or new inequalities created. Let $a_{new}$ be a newly generated train path for a train service $r$. The conflict graph $G = \left( \mathcal{A}^\prime, \mathcal{E} \right)$ is then expanded by this train path, i.e., $\mathcal{A}^\prime \leftarrow \mathcal{A}^\prime \cup \lbrace a_{new} \rbrace$. Through conflict detection \refSec{sec:conflict_detection}, conflicts between $a_{new}$ and $a \in \mathcal{A}^\prime \setminus \lbrace a_{new} \rbrace$ are identified, and corresponding edges $(a_{new},a)$ in $G$ are generated and added to the set of edges. In $G$ extended by $a_{new}$, maximal cliques $\mathcal{C}^\prime$ are now sought. Since only the maximal cliques involving $a_{new}$ are of interest, the search is performed in the subgraph of $G$ induced by $a_{new}$ denoted as $G \left( a_{new} \right) = \left( \mathcal{A}_{a_{new}}, \mathcal{E}_{a_{new}} \right)$. The set of nodes $\mathcal{A}_{a_{new}}$ consists of $a_{new}$ and and every $a \in \mathcal{A}^\prime$ having a conflict with a. For two nodes $a_1,a_2 \in \mathcal{A}_{a_{new}}$ there exists an edge $\left( a_1, a_2 \right)$ if $a_1,a_2$ have an conflict.
The detection of the set of the maximal cliques $\mathcal{C}^\prime$ in $G \left( a_{new} \right)$ is carried out using the algorithm by \citet{Tomita2017}, an extension of the algorithm by \citet{Bron1973}. Among these cliques $C^\prime \in \mathcal{C}^\prime$ containing $a_{new}$, a comparison is made to decide which clique $C \in \mathcal{C}$ needs to be expanded by $a_{new}$. For every $C^\prime$, a subclique $C = C^\prime \setminus \lbrace a_{new} \rbrace$ is searched, which is then extended by $a_{new}$, and the variable $x_{a_{new}}$ is added to the corresponding clique inequality \refEqu{equ:ip3} in the rRMP. If there is no such sub-clique $C$ for $C^\prime$ meeting the described conditions, then $C^\prime$ is a new maximal clique, and thus $\mathcal{C} \leftarrow \mathcal{C} \cup \lbrace C^\prime \rbrace$ follows along with the generation of a new clique inequality for $C^\prime$ in the rRMP, that is a row generation process.
\begin{figure}
	\centering
	\begin{subfigure}[t]{0.25\textwidth}
		\centering
		\input{img/clique_update_1}
	\end{subfigure}
	\begin{subfigure}[t]{0.35\textwidth}
		\centering
		\input{img/clique_update_2}
	\end{subfigure}
	\begin{subfigure}[t]{0.35\textwidth}
		\centering
		\input{img/clique_update_3}
	\end{subfigure}	
	\caption{Schematic representation of updating existing and detecting new maximal cliques. (a) Conflict graph without the new train path $a_{new}$. (b) Conflict graph with the new train path $a_{new}$. (c) Highlighted subgraph induced by $a_{new}$ with the updated (marked in red) and new (marked in green) maximal cliques.}
	\label{fig:clique_update}
\end{figure}
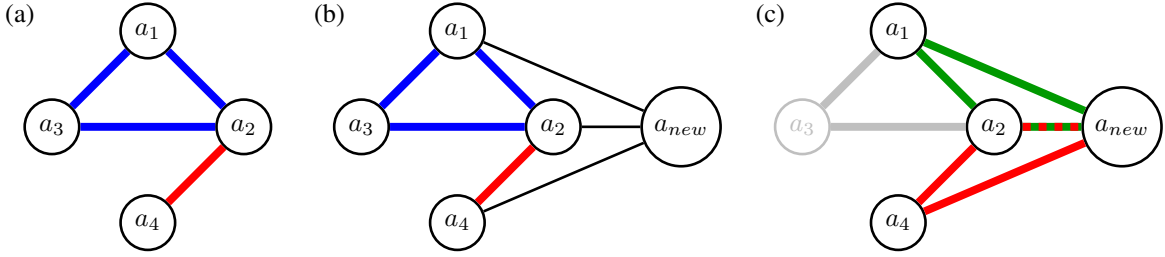

The process for updating the maximal cliques is depicted in \refFig{\ref{fig:clique_update}}. In \refFig[(a)]{\ref{fig:clique_update}}, the conflict graph $G$ shows the existing maximal cliques $C_1 = \lbrace a_1, a_2, a_3 \rbrace$ and $C_2 = \lbrace a_2, a_4 \rbrace$. A train path $a_{new}$ conflicts with the train paths $a_1$, $a_2$, and $a_4$ resulting in the extended conflict graph \refFigb[(b)]{\ref{fig:clique_update}}. In the subgraph of $G$ induced by $a_{new}$, the maximal cliques $C_1^\prime = \lbrace a_1, a_2, a_{new} \rbrace$ and $C_2^\prime = \lbrace a_2, a_4, a_{new} \rbrace$ exist \refFigb[(c)]{\ref{fig:clique_update}}. For $C_2$, only $C_2^\prime$ satisfies $\lvert C_2 \rvert = \lvert C_2^\prime \rvert - 1$, indicating that $C_2$ is a potentially updated clique. Moreover, since $\lvert C_2 \rvert = \lvert C_2^\prime \rvert \setminus \lbrace a_{new} \rbrace$ holds, $C_2$ will be updated. For $C_1^\prime$, neither $C_1$ nor $C_2$ meets these conditions, resulting in the addition of $C_1^\prime$ as a new clique.

With the new train path $a_{new}$, the updated and newly generated maximal cliques, the individual subproblems for the train services are also updated. For the new train path, the constraints \refEqu{equ:drive_up_conflict1} - \refEqu{equ:halting_conflict3} are generated in the subproblem \refEqu{equ:mip_subproblem} of every $r \in \mathcal{R}$. For each subproblem, for an updated clique $C$, the associated constraint \refEqu{equ:clique} is supplemented with the variable $z_{a_{new}}$, and for a newly generated maximal clique, such a constraint is added to the subproblem (row generation).

\section{Numerical Experiments} \label{sec:numerical_experiments}
\subsection{Scenario Description}
The dispatching area considered for the evaluation of the algorithm is located in the northern part of Bavaria, around Nuremberg (NN) \refFigb{\ref{fig:scenario}}.
Deutsche Bahn provides infrastructure data on the railway lines \citep{DB}. The Bavarian scenario was extracted from this database. The database includes informations on speed limits, block sections, junctions, single and double track lines \refFigb[(a),(b)]{\ref{fig:scenario_infra}}. From this, \profiles{} with minimum headway times can be constructed on the tracks. 
Due to the lack of information on the infrastructure within railway stations, we use a simplified transit matrix.
Arrival and departure nodes are connected via edges \refFigb[(c)]{\ref{fig:scenario_infra}}. For \profiles{} with geometrically crossing transit edges a minimum headway is generated within the preprocess. Halting positions on the transit edges were inserted to the best of our knowledge. 

Speed-profiles were constructed for regional and long distance trains for this scenario. Stopping patterns of the train services had been taken from a real-world timetable. Entry and exit nodes of the dispatching area  as well as halting points were defined as dispatching points, the \profiles{} were generated according to Section \ref{sec:speed_profiles}. 
Dispositive stops on open tracks are not considered in our scenario.
From a railway operations perspective, this leads to wasted capacity since trains are held at stations to give priority to other trains. Otherwise, trains could leave the station and let other trains pass at intersections, thereby not unnecessarily blocking platforms for extended periods. 
The focus of this work is on the methodological approach and its performance rather than on an exact optimisation of the real-world capacity.

In the dispatching area, we consider $n \in \lbrace 40, 50, 60\rbrace$ train services, including both long-distance and regional trains. For these train services with their scheduled stops, macroscopic travel paths are generated. 
Based on this, the \profiles{} are generated as described in Section \ref{sec:train_path_formulation}. 
\begin{figure}[ht!]
    \centering
    \scalebox{1.25}{\input{img/scenario}}
    \caption{Dispatching area in the region of northern Bavaria around Nuremberg.}
    \label{fig:scenario}
\end{figure}
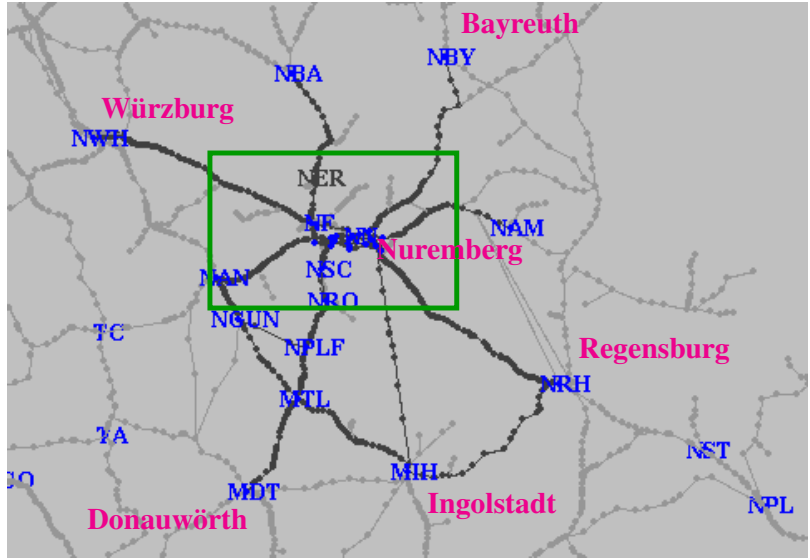 
\begin{figure}[ht!]
    \centering
    \scalebox{1.25}{\input{img/scenario_infra}}
	\caption{(a) A section of the operating point graph for the vicinity of Nuremberg Central Station (see \refFig{\ref{fig:scenario}}). (b) The infrastructure with block sections, single-track (red), and double-track (black and blue) lines based on the operating point graph in (a). (c) Connections for traversing Nuremberg Central Station.}

    \label{fig:scenario_infra}
\end{figure}
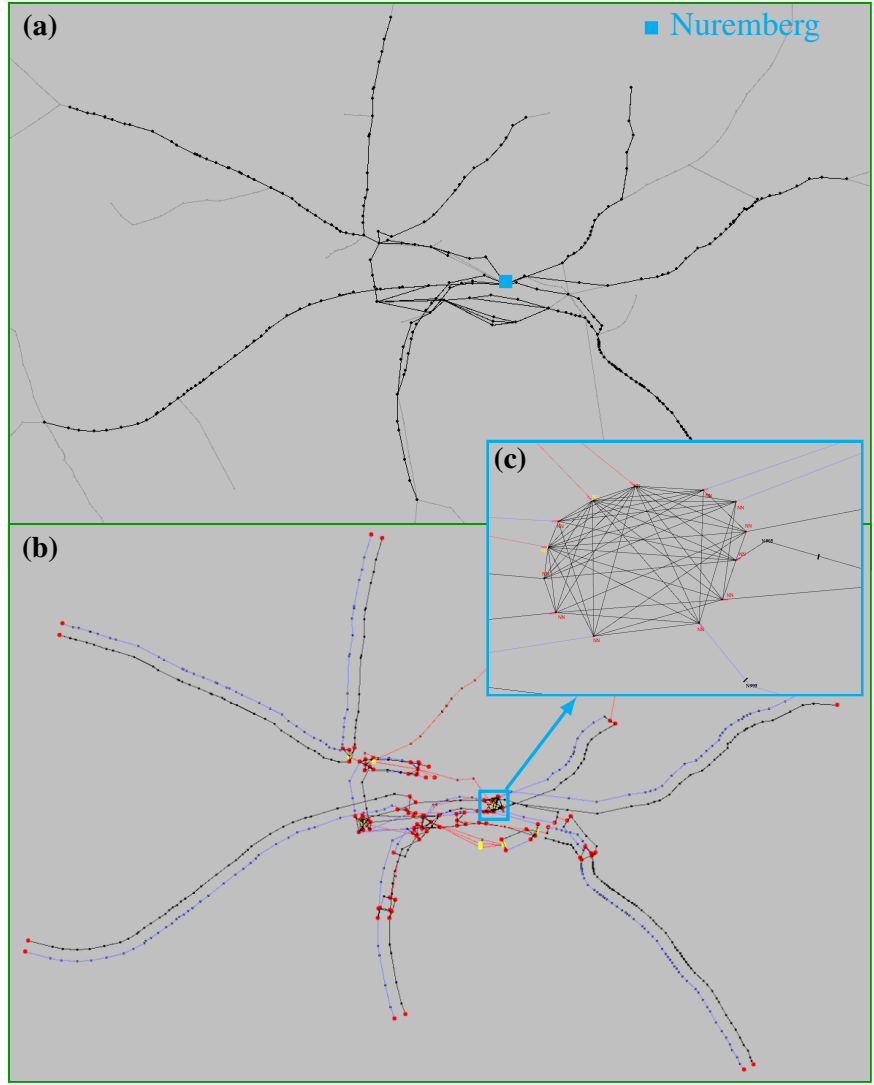 

In order to analyse the performance of the \ALGO{} with respect to different routing options, the full set $\mathcal{V} (r)$ generated according to Section \ref{sec:speed_profiles} is reduced by considering only the \profiles{} of at most $k$ shortest paths from entry to the exit points of the train service $r$. We call $k$ the reduction index. $k = \infty$ means, that we may use the full but finite set of \profiles{}. The path of the individual train services cut the dispatching area in different sizes, so that there are not always $k$ paths available.

\refTab{\ref{tab:results_complexity_routes}} displays the number of spacial routing options for the different instances\footnote{The actual number of available routing options within the reduced set of \profiles{} may be larger than the number of $k$ shortest path.}.
\begin{table}[ht!]
    \centering
    \caption{Total routing options for different number $n$ of train services when limiting the maximum routing options to $k$.}
    \begin{small}
		\input{table/results_complexity_routes_total}

    \end{small}
    \label{tab:results_complexity_routes}
\end{table}

\begin{figure}
	\centering
	\includegraphics[scale=0.75]{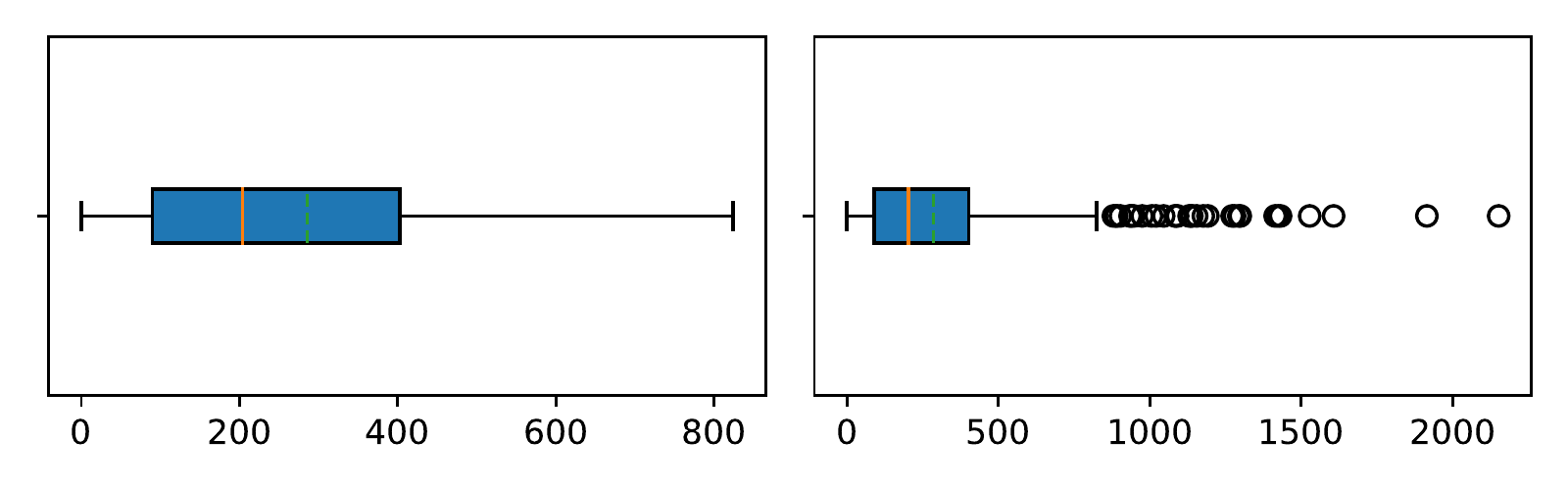}
	\caption{Boxplots of the initial disturbances for the scenarios.}
	\label{fig:boxplots}
\end{figure}

For the scenarios we consider train services with a regular, conflict-free timetable within a time horizon of approximately one hour.
The arrival times for these train services are subject to a disturbance $d$ which was modelled by the following probability function \citep{Schwanhaeusser1974}:
\begin{equation}
    P(D \leq d) =
    \begin{cases}
        0 & , \text{ if } d < 0 \\
        1 - \left( 1 - q \right) e^{-\lambda d} & , \text{ if } d \geq 0
    \end{cases}, \label{equ:delay}
\end{equation}
Positive delays follow an exponential distribution with the parameter $\lambda$ and $1 - q$ represents the proportion of disturbed train services. This approach gives us a railway operation situation for the TDP.
Disruption scenarios can range from (minor) disturbances to longer interruptions. From the boxplots in \refFig{\ref{fig:boxplots}} (sample size = 700), it can be seen that, for the scenarios, delays in the range of minor disturbances are considered. Due to the chosen distribution \refEqu{equ:delay} of the initial delays, the average is approximately 300 seconds in both cases. From the boxplots in the upper left and lower left in \refFig{\ref{fig:boxplots}}, it can be observed that 50\% of the values fall within the range of 100 to 400. For the scenarios, there are outliers above 1000 seconds.

In the Following, a scenario is encoded as NN-$n$-$k$. Here, NN designates the aforementioned dispatching area. For each scenario, 50 computations were performed with disturbances following Equation \refEqu{equ:delay} and fixed parameters $q = 0.8$ and $\lambda = \frac{1}{300} \frac{1}{\text{s}}$. All tests were conducted on an Intel Core i7-8550U CPU @ 1.80GHz with 16 GB RAM. The \ALGO{} is implemented in C++, and the IP, LP and MIP models are solved using the Gurobi solver (v9.0.1).

\subsection{Results}

\subsubsection{Complexity and Runtime Behaviour \ALGO}
In a first experiments we apply the \ALGO{} without any time limit in order to get the optimal solution. The number of train services in the dispatching area and the reduction index $k$ will be varied. We consider the scenarios NN-$n$-$k$ with $n \in \lbrace 40, 50, 60 \rbrace$ and $k \in \lbrace 100, 300, 500, 750, \infty \rbrace$. \refTab{\ref{tab:results_complexity_train_paths}} contains the average number of generated train paths. More train services in the scenarios lead to more train paths (variables $x_a$) being generated, whereas the reduction index does not have such a strong influence.
\begin{table}[ht!]
    \centering
    \caption{Average number of generated train paths (variables $x_a$) for scenarios with different $n$ and $k$.}
    \begin{small}
		\input{table/results_complexity_train_paths}

    \end{small}
    \label{tab:results_complexity_train_paths}
\end{table}

In \refTab{\ref{tab:results_complexity_cliques}}, the conflict situation for the individual scenarios is presented. As expected, this situation becomes more complex with a higher number of train services. The average number of cliques increases much more when the number of trains is increased, compared to an increase in routing options.
\begin{table}[ht!]
    \centering
    \caption{Conflict situation (average number of cliques / average clique size) for different maximum routing options and number of train services.}
    \begin{small}
		\input{table/results_complexity_cliques}

    \end{small}
    \label{tab:results_complexity_cliques}
\end{table}

For the scenarios NN-$n$-$k$ with $n \in \lbrace 40, 50, 60 \rbrace$ and $k \in \lbrace 100, 300, 500, 750, \infty \rbrace$, \refTab{\ref{tab:results_complexity_min_max_time}} shows the average, minimum and maximum computation time. It turns out that with a number of train services of 40, the average computation times are almost constant across different maximum routing options. Regarding the maximum computation time, it can be seen from \refTab{\ref{tab:results_complexity_min_max_time}} that there is an impact of the routing options for $n=50$ and $n=60$.  

As can be seen from \refTab{\ref{tab:results_complexity_routes}}, scenarios with 50 and 60 train services have almost identical total routing options. However, the computation times differ significantly \refTabb{\ref{tab:results_complexity_min_max_time}}. Hence, the number of train services has a greater influence on the runtime behavior. 

The higher computation times are associated with the conflict situation \refTabcf{\ref{tab:results_complexity_cliques} and \ref{tab:results_complexity_min_max_time}}. The increased computation times with more existing cliques result from the process of updating cliques. Here, a search for maximal cliques takes place, which is known to be a challenging problem. This can also explain the computation time for $n=60$ and $k=100$, as there were configurations that led to an \glqq unfavorable\grqq{} conflict situation. Here too, limiting the computation time is of crucial importance for real-world applications.
\begin{table}[ht!]
    \centering
    \caption{Computation times (average / minimum / maximum) [s] for different maximum routing options and number of train services.}
    \begin{small}
		\input{table/results_complexity_min_max_time}

    \end{small}
    \label{tab:results_complexity_min_max_time}
\end{table}

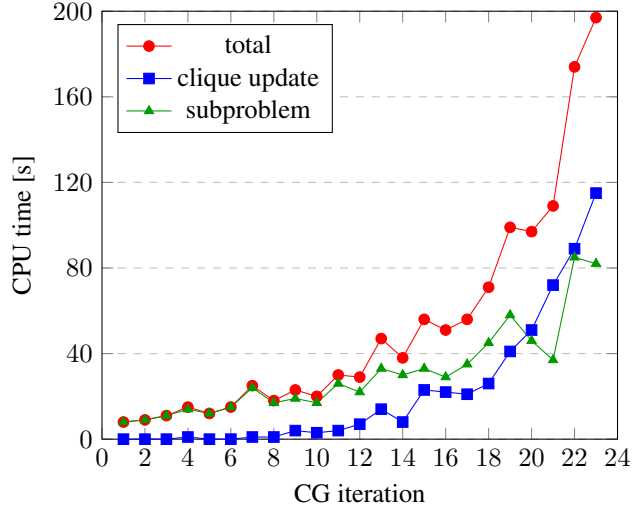
\begin{figure}[ht!]
	\centering
	\input{img/time_clique_detection}
	\caption{Computation time (CPU time) for the column generation iterations.}
	\label{fig:time_clique_detection}
\end{figure}

For the scenario NN-$60$-$\infty$ and the case with the maximum computation time of 1210 seconds, the computation times for each iteration of column generation are shown in \refFig{\ref{fig:time_clique_detection}}. The computation times grow exponentially. Hence, for real-world applications, it is advisable to terminate the column generation prematurely. The optimality of the solution is then not proven, but a high-quality solution may be obtained. To investigate this, the quality of the \ALGO{} will be evaluated in the following section.

\subsubsection{Quality of the \ALGO{} under Real-Time Requirements}

Further numerical computations aim to evaluate the quality of the algorithm under real-time conditions. For this purpose, two scenarios NN-$60$-$\infty$ are chosen, where the maximum computation time for the online capability of the algorithm is set to 60 seconds. In the first scenario, a minimum optimality gap (following only referred to as gap) of 0 is chosen, and in the second scenario, the gap is set to 0.1. Column generation suffers from the tailing-off effect \citep{Luebbecke2005}, so it is intended to terminate prematurely at a gap of 0.1. To assess the quality, the results of both scenarios will be compared with the scenario NN-$60$-$\infty$ with no restriction on the computation time and gap. The quality of the algorithm will be measured based on the delay quotient $\frac{d_{start}}{d_{end}}$. Here, $d_{start}$ is the total delay of the considered $n$ train services after the FCFS start heuristic, and $d_{end}$ is the delay after the termination of the column generation. 
\begin{table}[ht!]
    \centering
    \caption{Numerical results with and without computation time limits.}
    \begin{small}
		\input{table/results_quality}

    \end{small}
    \label{tab:results_quality}
\end{table}

\refTab{\ref{tab:results_quality}} presents the results for the described scenarios with different restriction to the minimum computation time and gap. Column 1 of \refTab{\ref{tab:results_quality}} contains the delay quotient. Column 2 indicates the average computation time (CPU time) of the algorithm (\ALGO{}). The optimality gap after the termination of the \ALGO{} (GAP) is listed in Column 3, and Column 4 contains the proportion of integer solutions after the termination of \ALGO{}. The penultimate column of \refTab{\ref{tab:results_quality}} displays the average number of cliques, and the last column lists the average number of iterations of column generation. Despite imposing constraints on the gap and computation time, the delay quotient demonstrates that the delay was reduced nearly equally well in all three scenarios, particularly for the scenarios with a zero gap. These restrictions significantly reduced the computation time, meeting real-time requirements. For the scenario with a gap of 0.1, the average computation time was further reduced by approximately 2.5 seconds, with only a minimal decrease in solution quality. \refTab{\ref{tab:results_quality}} further shows, that the \ALGO{} provided at least 92\% integer solutions. Again, the number of cliques corresponds with the computation time. In the scenario without constraints on the gap and computation time, an average of three to four additional iterations were performed. This indicates that in the final iterations of column generation, a large number of new cliques are generated, making the problem increasingly complex. This also makes it challenging to solve subproblems within an acceptable time \refFigcf{\ref{fig:time_clique_detection}}. Early branching is a common method for accelerating column generation \citep{Desaulniers1999}. However, after the pricing a high remaining optimality gap at the root node of the branching tree poses a greater risk that the branch-and-price algorithm may be less effective \citep{Maher2023}. Nevertheless, for the scenario with zero gap and a maximum computation time of 60 seconds, the gap is only 2\%, mitigating this risk, and a high-quality solution is already present at the root node.

\section{Conclusion} \label{sec:conclusion}

In this paper, a column generation approach for solving a path-oriented model for the Train Dispatching problem is presented. A major novelty is the modelling of conflicts as maximum cliques which have to be modified dynamicly during column generation. The resulting difficult problem of taking into account the shadow prices of cliques when calculating new train paths was formulated by a MIP and can be solved efficiently and quickly in practical use.
The strong formulation using path-based cliques provides good candidates for constructing integer solution after the column generation process. This statement supported by the observation that the relaxation very often provides integer solutions.

The separation between the calculation of speed profiles and the optimization of entire train paths allows a flexible use of infrastructure and driving dynamics calculations in different granulations and by different program systems.



The generation of \profiles{} based on real-world infrastructure is a topic for further research. With that, the \ALGO{} can and needs to be tested on more complex real-world instances. The shown computation times and the quality of the solutions promise good results even when applied to such instances. 

Another advantage of column generation results from its use in the context of a rolling horizon approach. With the \ALGO{}, the train paths are generated for the entire spatial observation space and can thus break out of the time horizon in time under certain circumstances. If the time horizon is now shifted in discrete time steps for the fixed spatial observation space (rolling horizon approach), the generated train paths retain their validity. Hence, the generated train paths can be used in subsequent optimization runs of the rolling horizon approach, where good lower bounds may already be available at the beginning of column generation, promising an advantage in terms of computation times and solution quality. 

Dispatching only takes place within the dispatching area. Possible effects beyond this area are not considered when applying dispatching measures. For the \ALGO{}, there is the possibility to formulate effects outside the dispatching area, such as connection or circulation conditions, in the form of cost functions and to integrate them into the objective function of the pricing problem.

\input{TD-GCM_bib_revision}
\end{document}

%% file: img/overview.tex
\tikzset{node/.style={circle, draw, line width=1pt, inner sep=0pt, text width=7mm, align=center}}
\tikzset{nodeF/.style={circle, fill=black, draw=black, line width=1pt, inner sep=0pt, text width=2mm, align=center}}
\tikzset{arc/.style={-latex, line width=1.0pt}}
\tikzset{sna/.style={->, decorate, decoration=snake, draw=black, line width=1.0pt}}
\tikzset{tra/.style={line width=1pt}}

\tikzset{dia/.style={diamond, aspect=5, draw,thick, minimum width=3cm, minimum height=0.5cm}}
\tikzset{rect/.style={draw, thick, minimum width=2cm, minimum height=1cm, align=center}}
\tikzset{rectRC/.style={draw, rounded corners, line width=2pt, minimum width=4.5cm, minimum height=1cm, align=center}}
\tikzset{trap/.style={trapezium, draw, thick, minimum width=3cm, minimum height=1cm, trapezium left angle=50, trapezium right angle=130, trapezium stretches}}

\tikzset{rectRCred/.style={draw, rounded corners, line width=2pt, minimum width=4.5cm, minimum height=1cm, align=center, draw=red, fill=red!10}}
\tikzset{rectRCblue/.style={draw, rounded corners, line width=2pt, minimum width=4.5cm, minimum height=1cm, align=center, draw=blue, fill=blue!10}}
\tikzset{rectRCgreen/.style={draw, rounded corners, line width=2pt, minimum width=2cm, minimum height=0.5cm, align=center, draw=lime green, fill=lime green!10}}
\tikzset{rectRCorange/.style={draw, rounded corners, line width=2pt, minimum width=4.5cm, minimum height=1cm, align=center, draw=orange, fill=orange!10}}
\tikzset{rectRClila/.style={draw, rounded corners, line width=2pt, minimum width=4.5cm, minimum height=1cm, align=center, draw=purple, fill=purple!10}}

\begin{tikzpicture}[scale=1]
	\def\width{7.0cm}
	\def\widthS{5.75cm}
	
	\node (prae) at (0,0) {\textbf{\textcolor{cyan}{preprocessing}} \refSec{sec:preprocessing}};
	\node[rect, below = 0.15cm of prae] (TP) {
		\begin{varwidth}{\width} 
			calculation and determination of:
			\begin{enumerate}[leftmargin=*]
				\item \profiles{} \refSec{sec:speed_profiles},
				\item conflict intervals \refSec{sec:conflict_detection} and
				\item potential holding conflicts \refSec{sec:conflict_detection}
			\end{enumerate}
		\end{varwidth}
	};
	
	\node[right = 4.25cm of prae] (state) {\textbf{\textcolor{red!60}{railway operating situation}}};
	\node[rect, below = 0.15cm of state] (STATE) {
		\begin{varwidth}{\widthS} 
			\begin{enumerate}[leftmargin=*]
				\item set of train services in dispatching area and horizon
				\item position of train services (spatial, temporal)
				\item scheduled timetable of train services 
			\end{enumerate}
		\end{varwidth}
	};
	
	\node[below = 1.5cm of TP] (online) {\textbf{\textcolor{limegreen}{\ALGO{} (online)}} \refSec{sec:PMC-CG}};
	\node[rect, below = 0.15cm of online] (IP) {
		\begin{varwidth}{\width} 
			\textbf{path-based IP} \refSec{sec:path_based_model} with
			\begin{enumerate}[leftmargin=*]
				\item path-based maximal clique inequalties
			\end{enumerate}
		\end{varwidth}
	};
	\node[rect, below = 1.5cm of IP] (CG) {
		\begin{varwidth}{\width} 
			\textbf{CG} \refSec{sec:column_generation}, where
			\begin{enumerate}[leftmargin=*]
				\item subproblems are formulated as MIPs
			\end{enumerate}
		\end{varwidth}
	};


	\begin{scope}[on background layer]
			\node[above = 0.25cm of online] (a1) {};
			\node[left = 0.25cm of IP] (l1) {};
			\node[below = 0.25cm of CG] (b1) {};
			\node[right = 0.25cm of IP] (r1) {};		
			\draw[limegreen, line width=1pt, fill=limegreen!5] (l1|-a1) rectangle (r1|-b1);
			\node[above = 0.25cm of prae] (a2) {};
			\node[left = 0.25cm of TP] (l2) {};
			\node[below = 0.25cm of TP] (b2) {};
			\node[right = 0.25cm of TP] (r2) {};
			\draw[cyan, line width=1pt, fill=cyan!5] (l2|-a2) rectangle (r2|-b2);	
			\node[above = 0.25cm of state] (a3) {};
			\node[left = 0.25cm of STATE] (l3) {};
			\node[below = 0.25cm of STATE] (b3) {};
			\node[right = 0.25cm of STATE] (r3) {};
			\draw[red!60, line width=1pt, fill=red!5] (l3|-a3) rectangle (r3|-b3);							
	\end{scope}

	\node[below = 0.75cm of IP] (h) {};

	\draw[arc] (IP) -- node[right] {solved by} (CG);
	\draw[arc] (prae|-b2) -- node[right] {used by} (online|-a1);
	\draw[arc] (state|-b3) to  [out=-90,in=0] node[below right] {used by} (r1|-h);
	\draw[arc, dashed] (l3|-STATE) -- node[above] {used} node[below] {by} (r2|-STATE);

\end{tikzpicture}

%% file: img/speed_profile_and_train_path.tex
\tikzset{dot/.style={circle, draw, fill, line width=1pt, inner sep=1pt}}
\tikzset{node/.style={circle, draw, line width=0.75pt, inner sep=0pt, text width=4.5mm, align=center}}
\tikzset{arc/.style={-latex, line width=0.75pt}}
\tikzset{line/.style={, line width=0.75pt}}

\def\blockL{1.5}
\def\platformL{2.5}
\def\crossingL{1}
\def\junctionL{1}
\def\d{1} 
\def\crossingOffa{0.125*\crossingL}
\def\crossingOffb{0.25*\crossingL}
\def\crossingOffc{0.75*\crossingL}
\def\crossingOffd{0.875*\crossingL}

\def\BlockB#1#2#3{
	\draw[](#1,#2+0.2)--(#1,#2-0.2);
	\draw[](#1,#2)--(#1+#3,#2);
	\draw[](#1+#3,#2+0.2)--(#1+#3,#2-0.2);
}

\def\Block#1#2#3#4{
	\draw[](#1,#2+0.2)--(#1,#2-0.2);
	\draw[](#1,#2)-- node[above] {\footnotesize{$b_{#4}$}} (#1+#3,#2);
	\draw[](#1+#3,#2+0.2)--(#1+#3,#2-0.2);
}

\def\BlockUp#1#2#3{
	\draw[](#1,#2+0.2)--(#1,#2-0.2);
	\draw[](#1,#2)-- node[above] {\footnotesize{$b_{#3}$}} (#1+\blockL,#2);
	\draw[](#1+\blockL,#2+0.2)--(#1+\blockL,#2-0.2);
}

\def\BlockDown#1#2#3{
	\draw[](#1,#2+0.2)--(#1,#2-0.2);
	\draw[](#1,#2)-- node[above] {\footnotesize{$b_{#3}$}} (#1+\blockL,#2);
	\draw[](#1+\blockL,#2+0.2)--(#1+\blockL,#2-0.2);
}

\def\JunctionLUp#1#2#3{
	\node (q) at (#1,#2) {};
	\node (z0) at (#1+\junctionL,#2) {};
	\node (z1) at (#1+\junctionL,#2+1) {};
	
	\draw[] (q.center) -- node[above, xshift=0.15cm] {\footnotesize{$b_{#3}$}} (z0.center);
	\draw[] (q.center) -- (z1.center);

	\draw[](#1,#2+0.2)--(#1,#2-0.2);
	\draw[](#1+\junctionL,#2+0.2)--(#1+\junctionL,#2-0.2);
	\draw[](#1+\junctionL,#2+1.2)--(#1+\junctionL,#2+0.8);
}

\def\JunctionLDown#1#2#3{
	\node (q) at (#1,#2) {};
	\node (z0) at (#1+\junctionL,#2) {};
	\node (z1) at (#1+\junctionL,#2-1) {};
	
	\draw[] (q.center) -- node[below, xshift=0.15cm] {\footnotesize{$b_{#3}$}} (z0.center);
	\draw[] (q.center) -- (z1.center);

	\draw[](#1,#2+0.2)--(#1,#2-0.2);
	\draw[](#1+\junctionL,#2+0.2)--(#1+\junctionL,#2-0.2);
	\draw[](#1+\junctionL,#2-1.2)--(#1+\junctionL,#2-0.8);
}

\def\JunctionL#1#2#3{
	\node (q) at (#1,#2) {};
	\node (z0) at (#1+\junctionL,#2) {};
	
	\draw[] (q.center) -- node[above, xshift=0.1cm] {\footnotesize{$b_{#3}$}} (z0.center);

	\draw[](#1,#2+0.2)--(#1,#2-0.2);
	\draw[](#1+\junctionL,#2+0.2)--(#1+\junctionL,#2-0.2);
}

\def\JunctionRUp#1#2#3{
	\node (q0) at (#1,#2) {};
	\node (q1) at (#1,#2+1) {};
	\node (z) at (#1+\junctionL,#2) {};
	
	\draw[] (q0.center) -- node[above, xshift=-0.15cm] {\footnotesize{$b_{#3}$}} (z.center);
	\draw[] (q1.center) -- (z.center);

	\draw[](#1,#2+0.2) -- (#1,#2-0.2);
	\draw[](#1,#2+1.2) -- (#1,#2+0.8);
	\draw[](#1+\junctionL,#2+0.2) -- (#1+\junctionL,#2-0.2);
}

\def\JunctionRDown#1#2#3{
	\node (q0) at (#1,#2) {};
	\node (q1) at (#1,#2-1) {};
	\node (z) at (#1+\junctionL,#2) {};
	
	\draw[] (q0.center) -- node[below, xshift=-0.15cm] {\footnotesize{$b_{#3}$}} (z.center);
	\draw[] (q1.center) -- (z.center);

	\draw[](#1,#2+0.2) -- (#1,#2-0.2);
	\draw[](#1,#2-1.2) -- (#1,#2-0.8);
	\draw[](#1+\junctionL,#2+0.2) -- (#1+\junctionL,#2-0.2);
}

\def\JunctionR#1#2#3{
	\node (q0) at (#1,#2) {};
	\node (z) at (#1+\junctionL,#2) {};
	
	\draw[] (q0.center) -- node[above, xshift=-0.1cm] {\footnotesize{$b_{#3}$}} (z.center);

	\draw[](#1,#2+0.2) -- (#1,#2-0.2);
	\draw[](#1+\junctionL,#2+0.2) -- (#1+\junctionL,#2-0.2);
}

\def\PlatformUp#1#2#3#4{
	\Block{#1}{#2}{\platformL}{#4};
	\Block{#1}{#2+1}{\platformL}{#3};
}

\def\PlatformDown#1#2#3#4{
	\Block{#1}{#2}{\platformL}{#3};
	\Block{#1}{#2-1}{\platformL}{#4};
}

\def\Platform#1#2#3{
	\Block{#1}{#2}{\platformL}{#3};
}

\def\Crossing#1#2#3{
	\BlockB{#1}{#2}{\crossingL};
	\Block{#1}{#2-\d}{\crossingL}{#3};
	\draw[] (#1+\crossingOffa,#2) -- (#1+\crossingOffb,#2-\d);
	\draw[] (#1+\crossingOffc,#2-\d) -- (#1+\crossingOffd,#2);
}

\begin{tikzpicture}[scale=0.5]
	
	\node[] at (0,4.5) {(a)};
	\node[node] (A) at(2.25+\platformL,3.5) {$S_1$};
	\node[node] (B) at(6.25+2*\platformL+2*\blockL,3.5) {$S_2$};
	\node[node] (C) at(10.25+3*\platformL+4*\blockL,3.5) {$S_3$};
	
	\draw[arc] (A) -- (B);
	\draw[arc] (B) -- (C);
	
	\Block{0.5}{0}{1}{0};
	\Block{0.5}{-\d}{1}{1};	
	\Crossing{\blockL}{0}{2};
	\JunctionLUp{\blockL+\crossingL}{0}{3};
	\JunctionL{\blockL+\crossingL}{-\d}{4};
	\PlatformUp{\junctionL+\blockL+\crossingL}{0}{5}{6};
	\Platform{\junctionL+\blockL+\crossingL}{-\d}{7}{8};
	\JunctionRUp{\junctionL+\blockL+\crossingL+\platformL}{0}{9};
	\JunctionR{\junctionL+\blockL+\crossingL+\platformL}{-\d}{10};
	\Crossing{2*\junctionL+\blockL+\crossingL+\platformL}{0}{11};
	\BlockUp{2*\junctionL+\blockL+2*\crossingL+\platformL}{0}{12};
	\BlockDown{2*\junctionL+\blockL+2*\crossingL+\platformL}{-\d}{13};
	\BlockUp{2*\junctionL+2*\blockL+2*\crossingL+\platformL}{0}{14};
	\BlockDown{2*\junctionL+2*\blockL+2*\crossingL+\platformL}{-\d}{15};
	\Crossing{2*\junctionL+3*\blockL+2*\crossingL+\platformL}{0}{16};
	\JunctionLUp{2*\junctionL+3*\blockL+3*\crossingL+\platformL}{0}{17};
	\JunctionLDown{2*\junctionL+3*\blockL+3*\crossingL+\platformL}{-\d}{18};
	\PlatformUp{3*\junctionL+3*\blockL+3*\crossingL+\platformL}{0}{19}{20};
	\PlatformDown{3*\junctionL+3*\blockL+3*\crossingL+\platformL}{-\d}{21}{22};
	\JunctionRUp{3*\junctionL+3*\blockL+3*\crossingL+2*\platformL}{0}{23};
	\JunctionRDown{3*\junctionL+3*\blockL+3*\crossingL+2*\platformL}{-\d}{24};
	\Crossing{4*\junctionL+3*\blockL+3*\crossingL+2*\platformL}{0}{25};
	\BlockUp{4*\junctionL+3*\blockL+4*\crossingL+2*\platformL}{0}{26};
	\BlockDown{4*\junctionL+3*\blockL+4*\crossingL+2*\platformL}{-\d}{27};
	\BlockUp{4*\junctionL+4*\blockL+4*\crossingL+2*\platformL}{0}{28};
	\BlockDown{4*\junctionL+4*\blockL+4*\crossingL+2*\platformL}{-\d}{29};
	\Crossing{4*\junctionL+5*\blockL+4*\crossingL+2*\platformL}{0}{30};
	\JunctionLUp{4*\junctionL+5*\blockL+5*\crossingL+2*\platformL}{0}{31};
	\JunctionL{4*\junctionL+5*\blockL+5*\crossingL+2*\platformL}{-\d}{32};
	\PlatformUp{5*\junctionL+5*\blockL+5*\crossingL+2*\platformL}{0}{33}{34};
	\Platform{5*\junctionL+5*\blockL+5*\crossingL+2*\platformL}{-\d}{35}{36};
	\JunctionRUp{5*\junctionL+5*\blockL+5*\crossingL+3*\platformL}{0}{37};
	\JunctionR{5*\junctionL+5*\blockL+5*\crossingL+3*\platformL}{-\d}{38};
	\Crossing{6*\junctionL+5*\blockL+5*\crossingL+3*\platformL}{0}{39};
	\Block{6*\junctionL+5*\blockL+6*\crossingL+3*\platformL}{0}{1}{40};
	\Block{6*\junctionL+5*\blockL+6*\crossingL+3*\platformL}{-\d}{1}{41};
	
	\draw[-latex,line width=1.5pt, blue] (\junctionL+\blockL+\crossingL+\platformL,-\d) -- (3*\junctionL+3*\blockL+3*\crossingL+2*\platformL,-\d);
	\draw[-latex,line width=1.5pt, cyan] (3*\junctionL+3*\blockL+3*\crossingL+2*\platformL,-\d) -- (5*\junctionL+5*\blockL+5*\crossingL+3*\platformL,-\d);

	\draw[-latex,line width=1.5pt, red, yshift=0.25cm] (3*\junctionL+3*\blockL+3*\crossingL+2*\platformL,-\d) -- (5*\junctionL+5*\blockL+5*\crossingL+3*\platformL,-\d);

	\draw[-latex,line width=1.5pt, green!60!black] (\junctionL+\blockL+\crossingL+\platformL,-\d-0.25) -- (2*\junctionL+3*\blockL+3*\crossingL+\platformL,-\d-0.25) -- (3*\junctionL+3*\blockL+3*\crossingL+\platformL,-\d-1.25) -- (3*\junctionL+3*\blockL+3*\crossingL+2*\platformL,-\d-1.25);
	\draw[-latex,line width=1.5pt, green] (3*\junctionL+3*\blockL+3*\crossingL+2*\platformL,-\d-1.25) -- (4*\junctionL+3*\blockL+3*\crossingL+2*\platformL,-\d-0.25) -- (5*\junctionL+5*\blockL+5*\crossingL+3*\platformL,-\d-0.25);

	\def\upshift{14}
	\node[] at (0,-17.25+\upshift) {(b)};
	\draw[arc] (-0.5,-18.5+\upshift) -- (29,-18.5+\upshift) node[right,xshift=-0.0cm] {\footnotesize{space}};
	\draw[arc] (0,-18+\upshift) -- (0,-33.0+\upshift) node[below] {\footnotesize{time}};
	
	\draw[arc, blue] (\junctionL+\blockL+\crossingL+\platformL,-19+\upshift) 
		to[out=-65, in=90, looseness=0.5, edge node={node [above] {$v_1$}}] 
		(3*\junctionL+3*\blockL+3*\crossingL+2*\platformL,-23+\upshift);
	\draw[arc, cyan] (3*\junctionL+3*\blockL+3*\crossingL+2*\platformL,-23+\upshift) 
		to[out=-65, in=90, looseness=0.5, edge node={node [above] {$v_2$}}] 
		(5*\junctionL+5*\blockL+5*\crossingL+3*\platformL,-27.5+\upshift);

	\draw[arc, blue] (\junctionL+\blockL+\crossingL+\platformL,-21+\upshift) 
		to[out=-65, in=90, looseness=0.5, edge node={node [above] {$v_1$}}] 
		(3*\junctionL+3*\blockL+3*\crossingL+2*\platformL,-25+\upshift);
	\draw[dashed, blue, line width=1pt] (3*\junctionL+3*\blockL+3*\crossingL+2*\platformL,-25+\upshift) -- (3*\junctionL+3*\blockL+3*\crossingL+2*\platformL,-27+\upshift);
	\draw[arc, red] (3*\junctionL+3*\blockL+3*\crossingL+2*\platformL,-27+\upshift) 
		to[out=-80, in=90, looseness=0.5, edge node={node [above] {$v_3$}}] 
		(5*\junctionL+5*\blockL+5*\crossingL+3*\platformL,-33+\upshift);

	\draw[green!60!black, line width=1pt] (\junctionL+\blockL+\crossingL+\platformL,-22+\upshift) 
		to[out=-65, in=160, looseness=0.5, edge node={node [above, xshift=0.92cm, yshift=-0.25cm] {$v_4$}}] 
		(3*\junctionL+3*\blockL+3*\crossingL+\platformL,-24.5+\upshift);
	\draw[arc, dashed, green!60!black] (3*\junctionL+3*\blockL+3*\crossingL+\platformL,-24.5+\upshift) 
		to[out=-10, in=90, looseness=0.75, edge node={node [above] {}}] 
		(3*\junctionL+3*\blockL+3*\crossingL+2*\platformL,-26+\upshift);
	\draw[arc, green] (3*\junctionL+3*\blockL+3*\crossingL+2*\platformL,-26+\upshift) 
		to[out=-65, in=90, looseness=0.5, edge node={node [above] {$v_5$}}] 
		(5*\junctionL+5*\blockL+5*\crossingL+3*\platformL,-30+\upshift);
		
	\draw[] (-0.2,-19+\upshift) -- node[right] {\footnotesize{$t_{v_1}$}} (0.2,-19+\upshift);
	\draw[] (-0.2,-21+\upshift) -- node[right] {\footnotesize{$t_{v_1}^\prime$}} (0.2,-21+\upshift);
	\draw[] (-0.2,-22+\upshift) -- node[right] {\footnotesize{$t_{v_4}$}} (0.2,-22+\upshift);
	\draw[] (-0.2,-23+\upshift) -- node[right] {\footnotesize{$t_{v_2}$}} (0.2,-23+\upshift);
	\draw[] (-0.2,-26+\upshift) -- node[right] {\footnotesize{$t_{v_5}$}} (0.2,-26+\upshift);
	\draw[] (-0.2,-27+\upshift) -- node[right] {\footnotesize{$t_{v_3}$}} (0.2,-27+\upshift);
	
	\draw[dotted, line width=0.75pt] (1.25,-19+\upshift) -- (\junctionL+\blockL+\crossingL+\platformL,-19+\upshift);
	\draw[dotted, line width=0.75pt] (1.25,-21+\upshift) -- (\junctionL+\blockL+\crossingL+\platformL,-21+\upshift);
	\draw[dotted, line width=0.75pt] (1.25,-22+\upshift) -- (\junctionL+\blockL+\crossingL+\platformL,-22+\upshift);
	\draw[dotted, line width=0.75pt] (1.25,-23+\upshift) -- (3*\junctionL+3*\blockL+3*\crossingL+2*\platformL,-23+\upshift);
	\draw[dotted, line width=0.75pt] (1.25,-26+\upshift) -- (3*\junctionL+3*\blockL+3*\crossingL+2*\platformL,-26+\upshift);
	\draw[dotted, line width=0.75pt] (1.25,-27+\upshift) -- (3*\junctionL+3*\blockL+3*\crossingL+2*\platformL,-27+\upshift);

	\draw[-latex] (-0.4,-19+\upshift) -- node[left, rotate=90, xshift=0.35cm, yshift=0.35cm] {$f_{v_1}$} (-0.4,-23+\upshift);
		
	\Block{0.5}{-18.5+\upshift}{1}{1};	
	\Block{1.5}{-18.5+\upshift}{1}{2};
	\Block{2.5}{-18.5+\upshift}{1}{4};
	\Block{3.5}{-18.5+\upshift}{\platformL}{7};
	\Block{3.5+\platformL}{-18.5+\upshift}{1}{10};
	\Block{4.5+\platformL}{-18.5+\upshift}{1}{11};
	\Block{5.5+\platformL}{-18.5+\upshift}{\blockL}{13};	
	\Block{5.5+\platformL+\blockL}{-18.5+\upshift}{\blockL}{15};
	\Block{5.5+\platformL+2*\blockL}{-18.5+\upshift}{1}{16};
	\Block{6.5+\platformL+2*\blockL}{-18.5+\upshift}{1}{18};
	\Block{7.5+\platformL+2*\blockL}{-18.5+\upshift}{\platformL}{21};
	\Block{7.5+2*\platformL+2*\blockL}{-18.5+\upshift}{1}{24};
	\Block{8.5+2*\platformL+2*\blockL}{-18.5+\upshift}{1}{25};
	\Block{9.5+2*\platformL+2*\blockL}{-18.5+\upshift}{\blockL}{27};
	\Block{9.5+2*\platformL+3*\blockL}{-18.5+\upshift}{\blockL}{29};
	\Block{9.5+2*\platformL+4*\blockL}{-18.5+\upshift}{1}{30};
	\Block{10.5+2*\platformL+4*\blockL}{-18.5+\upshift}{1}{32};
	\Block{11.5+2*\platformL+4*\blockL}{-18.5+\upshift}{\platformL}{35};
	\Block{11.5+3*\platformL+4*\blockL}{-18.5+\upshift}{1}{38};
	\Block{12.5+3*\platformL+4*\blockL}{-18.5+\upshift}{1}{39};
	\Block{13.5+3*\platformL+4*\blockL}{-18.5+\upshift}{1}{41};
\end{tikzpicture}

%% file: img/conflict_1.tex
\tikzset{node/.style={circle, draw, line width=1pt, inner sep=0pt, text width=7mm, align=center}}
\tikzset{dot/.style={circle, fill, draw, line width=1pt, inner sep=0pt, text width=3mm, align=center}}
\tikzset{arc/.style={-latex, line width=0.75pt}}

\def\Block#1#2#3{
	\draw[](#1,#2+0.2)--(#1,#2-0.2);
	\draw[](#1,#2)-- node[above] {$b_{#3}$} (#1+1.5,#2);
	\draw[](#1+1.5,#2+0.2)--(#1+1.5,#2-0.2);
}

\def\JunctionL#1#2#3{
	\node (q) at (#1,#2) {};
	\node (z0) at (#1+1,#2) {};
	\node (z1) at (#1+1,#2+1) {};
	
	\draw[] (q.center) -- node[above, xshift=0.15cm] {$b_{#3}$} (z0.center);
	\draw[] (q.center) -- (z1.center);

	\draw[](#1,#2+0.2)--(#1,#2-0.2);
	\draw[](#1+1,#2+0.2)--(#1+1,#2-0.2);
	\draw[](#1+1,#2+1.2)--(#1+1,#2+0.8);
}

\def\JunctionR#1#2#3{
	\node (q0) at (#1,#2) {};
	\node (q1) at (#1,#2+1) {};
	\node (z) at (#1+1,#2) {};
	
	\draw[] (q0.center) -- node[above, xshift=-0.15cm] {$b_{#3}$} (z.center);
	\draw[] (q1.center) -- (z.center);

	\draw[](#1,#2+0.2) -- (#1,#2-0.2);
	\draw[](#1,#2+1.2) -- (#1,#2+0.8);
	\draw[](#1+1,#2+0.2) -- (#1+1,#2-0.2);
}

\def\Platform#1#2#3#4{
	\Block{#1}{#2}{#3};
	\Block{#1}{#2+1}{#4};
}

\def\CurvedTrajectory#1#2#3#4#5{
	\draw [ultra thick,red] (-2,2) to[out=45,in=115] (1,1) to[out=-180+115,in=10] (-5,-3);
}

\def\refshift{1}
\def\lshift{0.825}
\def\vshift{1.24}
\def\wshift{1.25}
\def\wwshift{-0.51}

\begin{tikzpicture}[scale=1]
	
	
	\draw[->] (0,0) -- (14.0,0) node[right] {\footnotesize{$space$}};
	\draw[->] (0.25,0.25) -- (0.25,-4) node[below] {\footnotesize{$time$}};
	
	\draw[] (1.5,0.15) -- (1.5,-0.15);
	\node[] at (2,0.25) {\footnotesize{$b_1$}};
	\draw[] (2.5,0.15) -- (2.5,-0.15);
	\node[] at (3.25,0.25) {\footnotesize{$b_2$}};
	\draw[] (4.0,0.15) -- (4.0,-0.15);
	\node[] at (4.75,0.25) {\footnotesize{$b_3$}};
	\draw[] (5.5,0.15) -- (5.5,-0.15);
	\node[] at (6,0.25) {\footnotesize{$b_4$}};
	\draw[] (6.5,0.15) -- (6.5,-0.15);
	\node[] at (7.25,0.25) {\footnotesize{$b_5$}};
	\draw[] (8.0,0.15) -- (8.0,-0.15);	
	\node[] at (8.75,0.25) {\footnotesize{$b_6$}};
	\draw[] (9.5,0.15) -- (9.5,-0.15);
	\node[] at (10.25,0.25) {\footnotesize{$b_7$}};
	\draw[] (11,0.15) -- (11.0,-0.15);
	\node[] at (11.5,0.25) {\footnotesize{$b_8$}};
	\draw[] (12,0.15) -- (12,-0.15);
	\node[] at (12.75,0.25) {\footnotesize{$b_9$}};
	\draw[] (13.5,0.15) -- (13.5,-0.15);
	
	\draw[mygreen!50!white, line width=1.0pt] (1.5,-1.6-\refshift+\lshift+\wshift) rectangle ++(1.0,-0.8);
	\draw[mygreen!50!white, line width=1.0pt] (2.5,-1.9-\refshift+\lshift+\wshift) rectangle ++(1.5,-0.6);
	\draw[mygreen!50!white, line width=1.0pt] (4.0,-2.0-\refshift+\lshift+\wshift) rectangle ++(1.5,-0.9);
	
	\draw[arc, line width=1.5pt, mygreen] (1.5,-1.75-\refshift+\lshift+\wshift) to[out=-65, in=90, looseness=0.5, edge node={node [above] {$w$}}] (5.5,-2.75-\refshift+\lshift+\wshift);

	\draw[blue!50!white, line width=1.0pt] (1.5,-1.6-\refshift+\vshift) rectangle ++(1.0,-0.9);
	\draw[blue!50!white, line width=1.0pt] (2.5,-2.0-\refshift+\vshift) rectangle ++(1.5,-0.8);
	\draw[blue!50!white, line width=1.0pt] (4.0,-2.3-\refshift+\vshift) rectangle ++(1.5,-0.6);
	\draw[blue!50!white, line width=1.0pt] (5.5,-2.4-\refshift+\vshift) rectangle ++(1.0,-0.55);
	\draw[blue!50!white, line width=1.0pt] (6.5,-2.5-\refshift+\vshift) rectangle ++(1.5,-0.6);
	\draw[blue!50!white, line width=1.0pt] (8.0,-2.6-\refshift+\vshift) rectangle ++(1.5,-0.6);	
	\draw[blue!50!white, line width=1.0pt] (9.5,-2.8-\refshift+\vshift) rectangle ++(1.5,-0.5);
	\draw[blue!50!white, line width=1.0pt] (11,-2.9-\refshift+\vshift) rectangle ++(1.0,-0.5);
	\draw[blue!50!white, line width=1.0pt] (12,-3.0-\refshift+\vshift) rectangle ++(1.5,-1.1);

	\draw[arc, line width=1.5pt, blue] (1.5,-1.75-\refshift+\vshift) to[out=-65, in=90, looseness=0.25, edge node={node [above] {$v$}}] (13.5,-4-\refshift+\vshift);
	
	\draw[mygreen!50!white, line width=1.0pt] (1.5,-1.6-\refshift+\lshift+\wwshift) rectangle ++(1.0,-0.8);
	\draw[mygreen!50!white, line width=1.0pt] (2.5,-1.9-\refshift+\lshift+\wwshift) rectangle ++(1.5,-0.6);
	\draw[mygreen!50!white, line width=1.0pt] (4.0,-2.0-\refshift+\lshift+\wwshift) rectangle ++(1.5,-0.9);
	
	\draw[arc, line width=1.5pt, mygreen] (1.5,-1.75-\refshift+\lshift+\wwshift) to[out=-65, in=90, looseness=0.5, edge node={node [above] {$w$}}] (5.5,-2.75-\refshift+\lshift+\wwshift);

	\draw[] (0.35,-1.75-\refshift+\lshift+\wshift) -- (0.15,-1.75-\refshift+\lshift+\wshift) node[left, mygreen] {\footnotesize{$t_{w}$}};
	\draw[] (0.35,-1.75-\refshift+\vshift) -- (0.15,-1.75-\refshift+\vshift) node[left, blue] {\footnotesize{$t_{v} = 0$}};
	\draw[] (0.35,-1.75-\refshift+\lshift+\wwshift) -- (0.15,-1.75-\refshift+\lshift+\wwshift) node[left, mygreen] {\footnotesize{$t^\prime_{w}$}};

	\draw[] (0.35,-1.75-\refshift+\lshift+\wshift) -- (0.15,-1.75-\refshift+\lshift+\wshift);
	\draw[-latex] (0.25,-1.75-\refshift+\vshift) -- node[left] {\footnotesize{$-l \left( v,w \right)$}} (0.25,-1.75-\refshift+\lshift+\wshift);
	\draw[] (0.35,-1.75-\refshift+\vshift) -- (0.15,-1.75-\refshift+\vshift);
	\draw[-latex] (0.25,-1.75-\refshift+\vshift) -- node[left] {\footnotesize{$u \left( v,w \right)$}} (0.25,-1.75-\refshift+\lshift+\wwshift);
	\draw[] (0.35,-1.75-\refshift+\lshift+\wwshift) -- (0.15,-1.75-\refshift+\lshift+\wwshift);
	
\end{tikzpicture}

%% file: img/conflict_3.tex
\tikzset{node/.style={circle, draw, line width=1pt, inner sep=0pt, text width=7mm, align=center}}
\tikzset{dot/.style={circle, fill, draw, line width=1pt, inner sep=0pt, text width=3mm, align=center}}
\tikzset{arc/.style={-latex, line width=0.75pt}}

\def\Block#1#2#3{
	\draw[](#1,#2+0.2)--(#1,#2-0.2);
	\draw[](#1,#2)-- node[above] {$b_{#3}$} (#1+1.5,#2);
	\draw[](#1+1.5,#2+0.2)--(#1+1.5,#2-0.2);
}

\def\JunctionL#1#2#3{
	\node (q) at (#1,#2) {};
	\node (z0) at (#1+1,#2) {};
	\node (z1) at (#1+1,#2+1) {};
	
	\draw[] (q.center) -- node[above, xshift=0.15cm] {$b_{#3}$} (z0.center);
	\draw[] (q.center) -- (z1.center);

	\draw[](#1,#2+0.2)--(#1,#2-0.2);
	\draw[](#1+1,#2+0.2)--(#1+1,#2-0.2);
	\draw[](#1+1,#2+1.2)--(#1+1,#2+0.8);
}

\def\JunctionR#1#2#3{
	\node (q0) at (#1,#2) {};
	\node (q1) at (#1,#2+1) {};
	\node (z) at (#1+1,#2) {};
	
	\draw[] (q0.center) -- node[above, xshift=-0.15cm] {$b_{#3}$} (z.center);
	\draw[] (q1.center) -- (z.center);

	\draw[](#1,#2+0.2) -- (#1,#2-0.2);
	\draw[](#1,#2+1.2) -- (#1,#2+0.8);
	\draw[](#1+1,#2+0.2) -- (#1+1,#2-0.2);
}

\def\Platform#1#2#3#4{
	\Block{#1}{#2}{#3};
	\Block{#1}{#2+1}{#4};
}

\def\CurvedTrajectory#1#2#3#4#5{
	\draw [ultra thick,red] (-2,2) to[out=45,in=115] (1,1) to[out=-180+115,in=10] (-5,-3);
}

\def\refshift{1}
\def\vshift{1.17}
\def\wshift{2}
\def\wwshift{-0.21}

\begin{tikzpicture}[scale=1]

	\node[] at(-1,0) {$(a)$};
	
	\draw[->] (0,0) -- (14.0,0) node[right] {\footnotesize{$space$}};
	\draw[->] (0.25,0.25) -- (0.25,-5.25) node[below] {\footnotesize{$time$}};
	
	\draw[] (1.5,0.15) -- (1.5,-0.15);
	\node[] at (2,0.25) {\footnotesize{$b_1$}};
	\draw[] (2.5,0.15) -- (2.5,-0.15);
	\node[] at (3.25,0.25) {\footnotesize{$b_2$}};
	\draw[] (4.0,0.15) -- (4.0,-0.15);
	\node[] at (4.75,0.25) {\footnotesize{$b_3$}};
	\draw[] (5.5,0.15) -- (5.5,-0.15);
	\node[] at (6,0.25) {\footnotesize{$b_4$}};
	\draw[] (6.5,0.15) -- (6.5,-0.15);
	\node[] at (7.25,0.25) {\footnotesize{$b_5$}};
	\draw[] (8.0,0.15) -- (8.0,-0.15);	
	\node[] at (8.75,0.25) {\footnotesize{$b_6$}};
	\draw[] (9.5,0.15) -- (9.5,-0.15);
	\node[] at (10.25,0.25) {\footnotesize{$b_7$}};
	\draw[] (11,0.15) -- (11.0,-0.15);
	\node[] at (11.5,0.25) {\footnotesize{$b_8$}};
	\draw[] (12,0.15) -- (12,-0.15);
	\node[] at (12.75,0.25) {\footnotesize{$b_9$}};
	\draw[] (13.5,0.15) -- (13.5,-0.15);
	
	\draw[mygreen!50!white, line width=1.0pt] (1.5,-1.6-\refshift+\wshift) rectangle ++(1.0,-0.8);
	\draw[mygreen!50!white, line width=1.0pt] (2.5,-1.9-\refshift+\wshift) rectangle ++(1.5,-0.6);
	\draw[mygreen!50!white, line width=1.0pt] (4.0,-2.0-\refshift+\wshift) rectangle ++(1.5,-2.4);
	
	\draw[mygreen!50!white, line width=1.0pt] (5.5,-1.6-\refshift+\wwshift) rectangle ++(1.0,-0.85);
	\draw[mygreen!50!white, line width=1.0pt] (6.5,-2.1-\refshift+\wwshift) rectangle ++(1.5,-0.8);
	\draw[mygreen!50!white, line width=1.0pt] (8.0,-2.5-\refshift+\wwshift) rectangle ++(1.5,-1.1);
	
	\draw[arc, line width=1.5pt, mygreen] (1.5,-1.75-\refshift+\wshift) to[out=-65, in=90, looseness=0.5, edge node={node [above] {$w$}}] (5.5,-2.75-\refshift+\wshift);
	\draw[arc, line width=1.5pt, mygreen] (5.5,-1.75-\refshift+\wwshift) to[out=-65, in=90, looseness=0.5, edge node={node [above] {$w^\prime$}}] (9.5,-3.5-\refshift+\wwshift);

	\draw[blue!50!white, line width=1.0pt] (1.5,-1.6-\refshift+\vshift) rectangle ++(1.0,-0.9);
	\draw[blue!50!white, line width=1.0pt] (2.5,-2.0-\refshift+\vshift) rectangle ++(1.5,-0.8);
	\draw[blue!50!white, line width=1.0pt] (4.0,-2.3-\refshift+\vshift) rectangle ++(1.5,-0.6);
	\draw[blue!50!white, line width=1.0pt] (5.5,-2.4-\refshift+\vshift) rectangle ++(1.0,-0.55);
	\draw[blue!50!white, line width=1.0pt] (6.5,-2.5-\refshift+\vshift) rectangle ++(1.5,-0.6);
	\draw[blue!50!white, line width=1.0pt] (8.0,-2.6-\refshift+\vshift) rectangle ++(1.5,-0.6);	
	\draw[blue!50!white, line width=1.0pt] (9.5,-2.8-\refshift+\vshift) rectangle ++(1.5,-0.5);
	\draw[blue!50!white, line width=1.0pt] (11,-2.9-\refshift+\vshift) rectangle ++(1.0,-0.5);
	\draw[blue!50!white, line width=1.0pt] (12,-3.0-\refshift+\vshift) rectangle ++(1.5,-1.1);

	\draw[arc, line width=1.5pt, blue] (1.5,-1.75-\refshift+\vshift) to[out=-65, in=90, looseness=0.25, edge node={node [above] {$v$}}] (13.5,-4-\refshift+\vshift);

	\draw[line width=1.5pt, mygreen] (5.5,-2.75-\refshift+\wshift) -- node[xshift=-0.1cm] {\Huge{\textcolor{red}{\Lightning}}} (5.5,-1.75-\refshift+\wwshift);

	\draw[] (0.35,-1.75-\refshift+\wshift) -- (0.15,-1.75-\refshift+\wshift) node[left, mygreen] {\footnotesize{$t_{w}$}};
	\draw[] (0.35,-1.75-\refshift+\vshift) -- (0.15,-1.75-\refshift+\vshift) node[left, blue] {\footnotesize{$t_{v}$}};
	\draw[] (0.35,-1.75-\refshift+\wwshift) -- (0.15,-1.75-\refshift+\wwshift) node[left, mygreen] {\footnotesize{$t_{w^\prime}$}};
	
\end{tikzpicture}

%% file: img/conflict_5.tex
\tikzset{node/.style={circle, draw, line width=1pt, inner sep=0pt, text width=7mm, align=center}}
\tikzset{dot/.style={circle, fill, draw, line width=1pt, inner sep=0pt, text width=3mm, align=center}}
\tikzset{arc/.style={-latex, line width=0.75pt}}

\def\Block#1#2#3{
	\draw[](#1,#2+0.2)--(#1,#2-0.2);
	\draw[](#1,#2)-- node[above] {$b_{#3}$} (#1+1.5,#2);
	\draw[](#1+1.5,#2+0.2)--(#1+1.5,#2-0.2);
}

\def\JunctionL#1#2#3{
	\node (q) at (#1,#2) {};
	\node (z0) at (#1+1,#2) {};
	\node (z1) at (#1+1,#2+1) {};
	
	\draw[] (q.center) -- node[above, xshift=0.15cm] {$b_{#3}$} (z0.center);
	\draw[] (q.center) -- (z1.center);

	\draw[](#1,#2+0.2)--(#1,#2-0.2);
	\draw[](#1+1,#2+0.2)--(#1+1,#2-0.2);
	\draw[](#1+1,#2+1.2)--(#1+1,#2+0.8);
}

\def\JunctionR#1#2#3{
	\node (q0) at (#1,#2) {};
	\node (q1) at (#1,#2+1) {};
	\node (z) at (#1+1,#2) {};
	
	\draw[] (q0.center) -- node[above, xshift=-0.15cm] {$b_{#3}$} (z.center);
	\draw[] (q1.center) -- (z.center);

	\draw[](#1,#2+0.2) -- (#1,#2-0.2);
	\draw[](#1,#2+1.2) -- (#1,#2+0.8);
	\draw[](#1+1,#2+0.2) -- (#1+1,#2-0.2);
}

\def\Platform#1#2#3#4{
	\Block{#1}{#2}{#3};
	\Block{#1}{#2+1}{#4};
}

\def\CurvedTrajectory#1#2#3#4#5{
	\draw [ultra thick,red] (-2,2) to[out=45,in=115] (1,1) to[out=-180+115,in=10] (-5,-3);
}

\def\refshift{1}
\def\vshift{1.17}
\def\vvshift{-0.5}
\def\wshift{2}
\def\wwshift{-1.25}

\begin{tikzpicture}[scale=1]

	\node[] at(-1,0) {$(b)$};
	
	\draw[->] (0,0) -- (14.0,0) node[right] {\footnotesize{$space$}};
	\draw[->] (0.25,0.25) -- (0.25,-6.25) node[below] {\footnotesize{$time$}};
	
	\draw[] (1.5,0.15) -- (1.5,-0.15);
	\node[] at (2,0.25) {\footnotesize{$b_1$}};
	\draw[] (2.5,0.15) -- (2.5,-0.15);
	\node[] at (3.25,0.25) {\footnotesize{$b_2$}};
	\draw[] (4.0,0.15) -- (4.0,-0.15);
	\node[] at (4.75,0.25) {\footnotesize{$b_3$}};
	\draw[] (5.5,0.15) -- (5.5,-0.15);
	\node[] at (6,0.25) {\footnotesize{$b_4$}};
	\draw[] (6.5,0.15) -- (6.5,-0.15);
	\node[] at (7.25,0.25) {\footnotesize{$b_5$}};
	\draw[] (8.0,0.15) -- (8.0,-0.15);	
	\node[] at (8.75,0.25) {\footnotesize{$b_6$}};
	\draw[] (9.5,0.15) -- (9.5,-0.15);
	\node[] at (10.25,0.25) {\footnotesize{$b_7$}};
	\draw[] (11,0.15) -- (11.0,-0.15);
	\node[] at (11.5,0.25) {\footnotesize{$b_8$}};
	\draw[] (12,0.15) -- (12,-0.15);
	\node[] at (12.75,0.25) {\footnotesize{$b_9$}};
	\draw[] (13.5,0.15) -- (13.5,-0.15);
	
	\draw[mygreen!50!white, line width=1.0pt] (1.5,-1.6-\refshift+\wshift) rectangle ++(1.0,-0.8);
	\draw[mygreen!50!white, line width=1.0pt] (2.5,-1.9-\refshift+\wshift) rectangle ++(1.5,-0.6);
	\draw[mygreen!50!white, line width=1.0pt] (4.0,-2.0-\refshift+\wshift) rectangle ++(1.5,-3.4);
	
	\draw[mygreen!50!white, line width=1.0pt] (5.5,-1.6-\refshift+\wwshift) rectangle ++(1.0,-0.85);
	\draw[mygreen!50!white, line width=1.0pt] (6.5,-2.1-\refshift+\wwshift) rectangle ++(1.5,-0.8);
	\draw[mygreen!50!white, line width=1.0pt] (8.0,-2.5-\refshift+\wwshift) rectangle ++(1.5,-1.1);
	
	\draw[arc, line width=1.5pt, mygreen] (1.5,-1.75-\refshift+\wshift) to[out=-65, in=90, looseness=0.5, edge node={node [above] {$w$}}] (5.5,-2.75-\refshift+\wshift);
	\draw[arc, line width=1.5pt, mygreen] (5.5,-1.75-\refshift+\wwshift) to[out=-65, in=90, looseness=0.5, edge node={node [above] {$w^\prime$}}] (9.5,-3.5-\refshift+\wwshift);

	\draw[blue!50!white, line width=1.0pt] (1.5,-1.6-\refshift+\vshift) rectangle ++(1.0,-0.8);
	\draw[blue!50!white, line width=1.0pt] (2.5,-1.9-\refshift+\vshift) rectangle ++(1.5,-0.6);
	\draw[blue!50!white, line width=1.0pt] (4.0,-2.0-\refshift+\vshift) rectangle ++(1.5,-1.8);
	
	\draw[blue!50!white, line width=1.0pt] (5.5,-1.6-\refshift+\vvshift) rectangle ++(1.0,-0.85);
	\draw[blue!50!white, line width=1.0pt] (6.5,-2.1-\refshift+\vvshift) rectangle ++(1.5,-0.8);
	\draw[blue!50!white, line width=1.0pt] (8.0,-2.5-\refshift+\vvshift) rectangle ++(1.5,-1.1);	

	\draw[arc, line width=1.5pt, blue] (1.5,-1.75-\refshift+\vshift) to[out=-65, in=90, looseness=0.5, edge node={node [above] {$v$}}] (5.5,-2.75-\refshift+\vshift);
	\draw[arc, line width=1.5pt, blue] (5.5,-1.75-\refshift+\vvshift) to[out=-65, in=90, looseness=0.5, edge node={node [above] {$v^\prime$}}] (9.5,-3.5-\refshift+\vvshift);

	\draw[line width=1.5pt, mygreen] (5.5,-2.75-\refshift+\wshift) -- node[xshift=-0.1cm] {\Huge{\textcolor{red}{\Lightning}}} (5.5,-1.75-\refshift+\wwshift);

	\draw[] (0.35,-1.75-\refshift+\wshift) -- (0.15,-1.75-\refshift+\wshift) node[left, mygreen] {\footnotesize{$t_{w}$}};
	\draw[] (0.35,-1.75-\refshift+\wwshift) -- (0.15,-1.75-\refshift+\wwshift) node[left, mygreen] {\footnotesize{$t_{w^\prime}$}};
	\draw[] (0.35,-1.75-\refshift+\vshift) -- (0.15,-1.75-\refshift+\vshift) node[left, blue] {\footnotesize{$t_{v}$}};
	\draw[] (0.35,-1.75-\refshift+\vvshift) -- (0.15,-1.75-\refshift+\vvshift) node[left, blue] {\footnotesize{$t_{v^\prime}$}};

\end{tikzpicture}

%% file: img/algorithm_overview.tex
\tikzset{node/.style={circle, draw, line width=1pt, inner sep=0pt, text width=7mm, align=center}}
\tikzset{nodeF/.style={circle, fill=black, draw=black, line width=1pt, inner sep=0pt, text width=2mm, align=center}}
\tikzset{arc/.style={-latex, line width=1.0pt}}
\tikzset{sna/.style={->, decorate, decoration=snake, draw=black, line width=1.0pt}}
\tikzset{tra/.style={line width=1pt}}

\tikzset{dia/.style={diamond, aspect=3, draw,thick, minimum width=2cm, minimum height=0.5cm}}
\tikzset{rect/.style={draw, thick, minimum width=2cm, minimum height=1cm, align=center}}
\tikzset{rectRC/.style={draw, rounded corners=14pt, thick, minimum width=2cm, minimum height=1cm, align=center}}
\tikzset{trap/.style={trapezium, draw, thick, minimum width=3cm, minimum height=1cm, trapezium left angle=50, trapezium right angle=130, trapezium stretches}}

\tikzset{rectRCred/.style={draw, rounded corners, line width=2pt, minimum width=4.5cm, minimum height=1cm, align=center, draw=red, fill=red!10}}
\tikzset{rectRCblue/.style={draw, rounded corners, line width=2pt, minimum width=4.5cm, minimum height=1cm, align=center, draw=blue, fill=blue!10}}
\tikzset{rectRCgreen/.style={draw, rounded corners, line width=2pt, minimum width=2cm, minimum height=0.5cm, align=center, draw=lime green, fill=lime green!10}}
\tikzset{rectRCorange/.style={draw, rounded corners, line width=2pt, minimum width=4.5cm, minimum height=1cm, align=center, draw=orange, fill=orange!10}}
\tikzset{rectRClila/.style={draw, rounded corners, line width=2pt, minimum width=4.5cm, minimum height=1cm, align=center, draw=purple, fill=purple!10}}

\begin{tikzpicture}[scale=0.75]
	\def\sep{0.75cm}
	\def\sepR{7cm}
	\def\sepRr{6.5cm}
	\def\sepL{8cm}
	\def\sepLl{12cm}
	
	\node[trap] (start_solution) at (0,0) {start solution};
	\node[rectRC] (start) at([shift=({-\sepL,0cm})]start_solution) {start};
	\node[rect, below = \sep of start_solution] (RMP) {solve the rRMP};
	\node[trap, below = \sep of RMP] (dual_solution) {
		\begin{varwidth}{\linewidth}
			\begin{center}
				solution of the rRMP and \\
				$\alpha_r$ for all $r \in \mathcal{R} $, $\beta_C$ for all $C \in \mathcal{C}$
			\end{center}
		\end{varwidth}	
	};
	\node[rect, below = \sep of dual_solution] (subproblem) {solve \refEqu{equ:mip_subproblem} for all $r \in \mathcal{R}$};
	\node[dia, below = \sep of subproblem] (conditions) {
		\begin{varwidth}{\linewidth}
			\begin{center}
				optimality gap / \\
				time limit reached? 
			\end{center}
		\end{varwidth}	
	};
	\node[dia, below = \sep of conditions] (negative_reduced_cost) {$a_{new}$ with $\rho_{a_{new}} < 0$?};
	\node[dia, below = \sep of negative_reduced_cost] (integer_fractional) {integer solution?};
	\node[rect] (get_integer) at([shift=({-\sepL,0cm})]integer_fractional) {solve the RMP};
	\node[rectRC, below = \sep of integer_fractional] (done_one) {finished};
	\node[dia] (integer_fractional_two) at([shift=({-\sepL,0cm})]conditions) {integer solution?};	
	
	\node[rect] (update_RMP) at([shift=({\sepRr,0cm})]subproblem) {
		\begin{varwidth}{\linewidth}
			\begin{center}
				update the rRMP \\
				and subproblems \\
				by row and \\
				column generation
			\end{center}
		\end{varwidth}				
	};
	\node[] (h1) at ([shift=({\sepRr,0cm})]RMP) {};
	\node[] (h2) at([shift=({\sepRr,0cm})]negative_reduced_cost) {};
	\node[] (h3) at([shift=({-\sepL,0cm})]done_one) {};
	\node[] (h4) at([shift=({-\sepLl,0cm})]conditions) {};
	\node[] (h5) at([shift=({-\sepLl,0cm})]done_one) {};

	\draw[arc] (start) -- (start_solution);	
	\draw[arc] (start_solution) -- (RMP);	
	\draw[arc] (RMP) -- (dual_solution);
	\draw[arc] (dual_solution) -- (subproblem);
	\draw[arc] (subproblem) -- (conditions);
	\draw[arc] (conditions) -- node[right] {no} (negative_reduced_cost);
	\draw[arc] (negative_reduced_cost) -- node[right] {no} (integer_fractional);
	\draw[arc] (integer_fractional) -- node[above] {no} (get_integer);
	\draw[arc] (conditions) -- node[above] {yes} (integer_fractional_two);
	\draw[arc] (integer_fractional_two) -- node[right] {no} (get_integer);
	\draw[arc] (integer_fractional) -- node[right] {yes} (done_one);
	
	\draw[arc] (get_integer) -- (h3.center) -- (done_one);
	\draw[] (integer_fractional_two) -- node[above] {yes} (h4.center);
	\draw[arc] (h4.center) -- (h5.center) -- (done_one);	
	\draw[] (negative_reduced_cost) -- node[above] {yes} (h2.center);
	\draw[arc] (h2.center) -- (update_RMP); 
	\draw[] (update_RMP) -- (h1.center);
	\draw[arc] (h1.center) -- (RMP);
\end{tikzpicture}

%% file: img/clique_update_1.tex
\definecolor{lime green}{RGB}{50,205,50}

\tikzset{node/.style={circle, draw=black, line width=1pt}}
\tikzset{fnode/.style={circle, draw=black, fill=black, line width=1pt, inner sep=2pt}}
\tikzset{rec/.style={rectangle, draw=black, line width=1pt}}
\tikzset{elip/.style={ellipse, minimum height=0.75cm, minimum width=1.5cm,draw}}
\tikzset{arc/.style={->, line width=1pt}}
\tikzset{disturbed/.style={dashed, lime green, line width=1pt}}
\tikzset{>=latex}

\begin{tikzpicture}[scale=0.85]

	\node[] at(-2,0.25) {(a)};

	\node[node] (a1) at (0,0) {$a_1$};
	\node[node] (a2) at (1.5,-1.5) {$a_2$};
	\node[node] (a3) at (-1.5,-1.5) {$a_3$};
	\node[node] (a4) at (0,-3) {$a_4$};
	
	\draw[line width = 3pt, blue] (a1) -- (a2);
	\draw[line width = 3pt, blue] (a1) -- (a3);
	\draw[line width = 3pt, blue] (a2) -- (a3);
	\draw[line width = 3pt, red] (a2) -- (a4);

\end{tikzpicture}

%% file: img/clique_update_2.tex
\definecolor{lime green}{RGB}{50,205,50}

\tikzset{node/.style={circle, draw=black, line width=1pt}}
\tikzset{fnode/.style={circle, draw=black, fill=black, line width=1pt, inner sep=2pt}}
\tikzset{rec/.style={rectangle, draw=black, line width=1pt}}
\tikzset{elip/.style={ellipse, minimum height=0.75cm, minimum width=1.5cm,draw}}
\tikzset{arc/.style={->, line width=1pt}}
\tikzset{disturbed/.style={dashed, lime green, line width=1pt}}
\tikzset{>=latex}

\begin{tikzpicture}[scale=0.85]

	\node[] at(-2,0.25) {(b)};

	\node[node] (a1) at (0,0) {$a_1$};
	\node[node] (a2) at (1.5,-1.5) {$a_2$};
	\node[node] (a3) at (-1.5,-1.5) {$a_3$};
	\node[node] (a4) at (0,-3) {$a_4$};
	
	\draw[line width = 3pt, blue] (a1) -- (a2);
	\draw[line width = 3pt, blue] (a1) -- (a3);
	\draw[line width = 3pt, blue] (a2) -- (a3);
	\draw[line width = 3pt, red] (a2) -- (a4);
	
	\node[node] (anew) at (3.5,-1.5) {$a_{new}$};
	
	\draw[line width = 1pt] (a1) -- (anew);
	\draw[line width = 1pt] (a2) -- (anew);
	\draw[line width = 1pt] (a4) -- (anew);

\end{tikzpicture}

%% file: img/clique_update_3.tex
\definecolor{lime green}{RGB}{50,205,50}

\tikzset{node/.style={circle, draw=black, line width=1pt}}
\tikzset{fnode/.style={circle, draw=black, fill=black, line width=1pt, inner sep=2pt}}
\tikzset{rec/.style={rectangle, draw=black, line width=1pt}}
\tikzset{elip/.style={ellipse, minimum height=0.75cm, minimum width=1.5cm,draw}}
\tikzset{arc/.style={->, line width=1pt}}
\tikzset{disturbed/.style={dashed, lime green, line width=1pt}}
\tikzset{>=latex}

\begin{tikzpicture}[scale=0.85]

	\node[] at(-2,0.25) {(c)};

	\node[node] (a1) at (0,0) {$a_1$};
	\node[node] (a2) at (1.5,-1.5) {$a_2$};
	\node[node, gray!50] (a3) at (-1.5,-1.5) {$a_3$};
	\node[node] (a4) at (0,-3) {$a_4$};
	
	\draw[line width = 3pt, mygreen] (a1) -- (a2);
	\draw[line width = 3pt, gray!50] (a1) -- (a3);
	\draw[line width = 3pt, gray!50] (a2) -- (a3);
	\draw[line width = 3pt, red] (a2) -- (a4);
	
	\node[node] (anew) at (3.5,-1.5) {$a_{new}$};
	
	\draw[line width = 3pt, mygreen] (a1) -- (anew);
	\draw[line width = 3pt, mygreen] (a2) -- (anew);
	\draw[line width = 3pt, red, dashed] (a2) -- (anew);
	\draw[line width = 3pt, red] (a4) -- (anew);

\end{tikzpicture}

%% file: img/scenario.tex
\tikzset{rectEC/.style={fill=blue!10, draw, thick, minimum width=1cm, minimum height=0.5cm}}
\tikzset{rectRE/.style={fill=green!10, draw, thick, minimum width=1cm, minimum height=0.5cm}}
\tikzset{rectICE/.style={fill=red!10, draw, thick, minimum width=1cm, minimum height=0.5cm}}
\tikzset{rectUm/.style={fill=yellow!10, draw, thick, minimum width=1cm, minimum height=0.5cm}}

\begin{tikzpicture}[scale=1.50]
	\node[inner sep=0pt, draw=none, fill=none] (image) at (0,0) 
		{\includegraphics[trim={0 0 0 0},clip,scale=0.5]{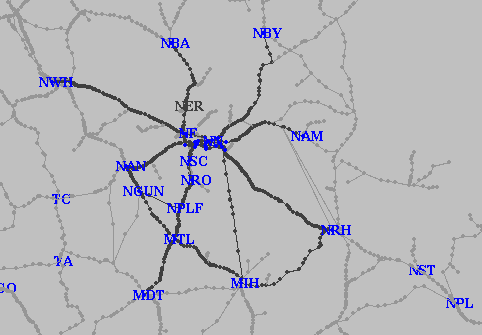}};

	\draw[draw=mygreen, line width=1.5pt, fill=none] (-1.4,-0.2) rectangle (0.35,0.9);
		
	\node[magenta] at(0.3,0.2) {\footnotesize{\textbf{Nuremberg}}};
	\node[magenta] at(-1.7,1.2) {\footnotesize{\textbf{Würzburg}}};
	\node[magenta] at(0.8,1.8) {\footnotesize{\textbf{Bayreuth}}};
	\node[magenta] at(-1.7,-1.7) {\footnotesize{\textbf{Donauwörth}}};
	\node[magenta] at(0.6,-1.6) {\footnotesize{\textbf{Ingolstadt}}};
	\node[magenta] at(1.75,-0.5) {\footnotesize{\textbf{Regensburg}}};

\end{tikzpicture}

%% file: img/scenario_infra.tex
\begin{tikzpicture}[scale=1.0]
	\node[draw=mygreen, line width=1.5pt, inner sep=0pt, fill=none] (bst) at (0,0) 
		{\includegraphics[trim={0 0 0 0},clip,scale=0.22]{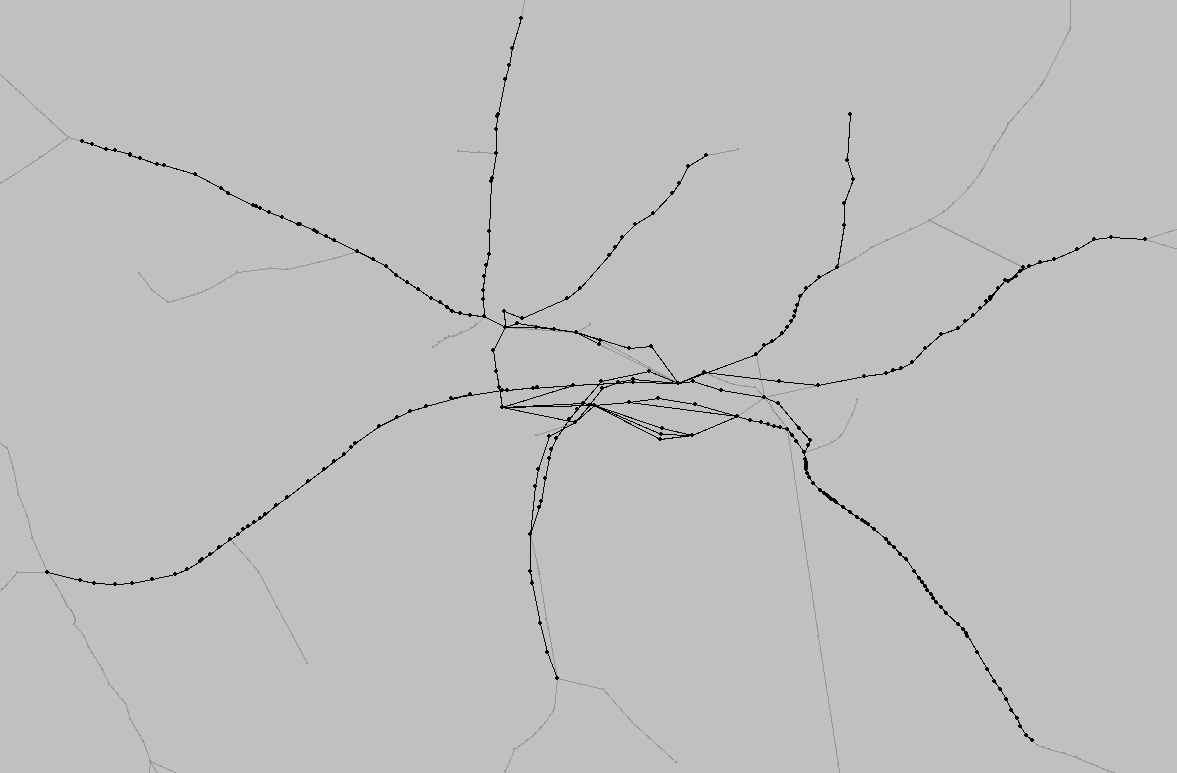}};
	\node[draw=mygreen, line width=1.5pt, inner sep=0pt, fill=none] (infra) at (0,-5.5) 
		{\includegraphics[trim={0 0 0 0},clip,scale=0.22]{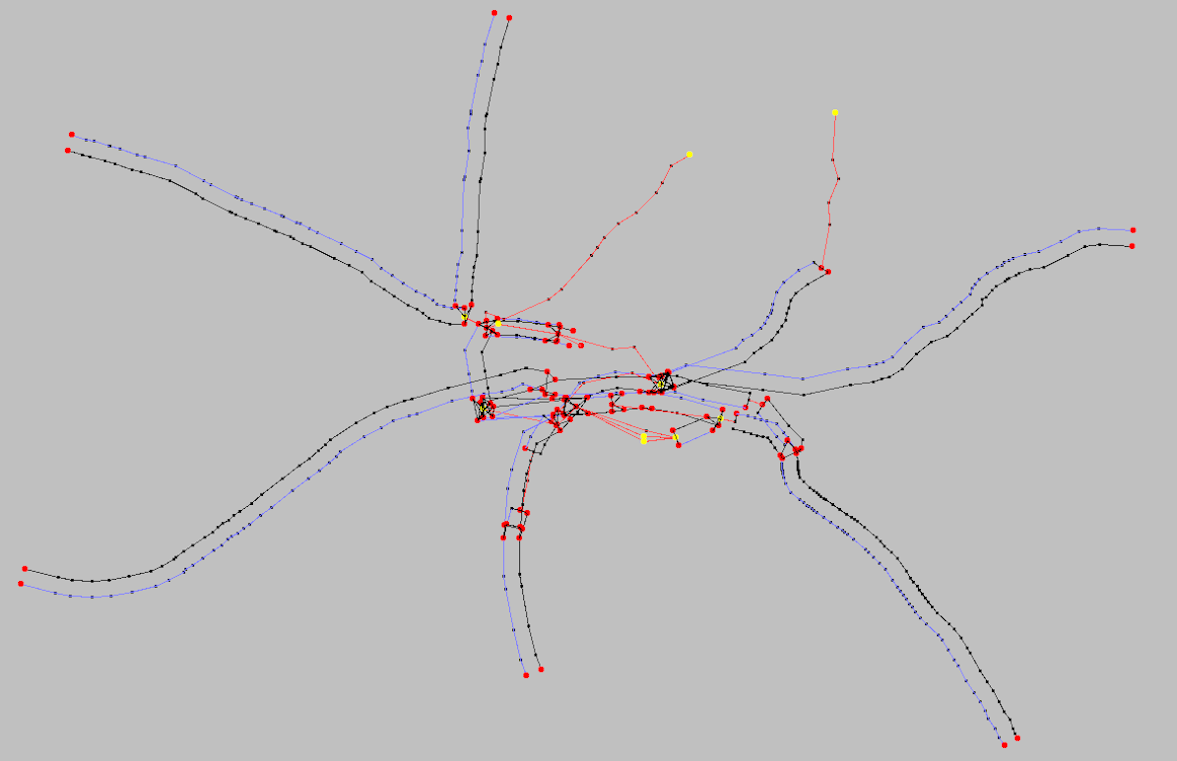}};
	\node[draw=cyan, line width=2pt, inner sep=0pt, fill=none] (NNzoom) at (2.5,-3) 
		{\includegraphics[trim={0 0 0 0},clip,scale=0.10]{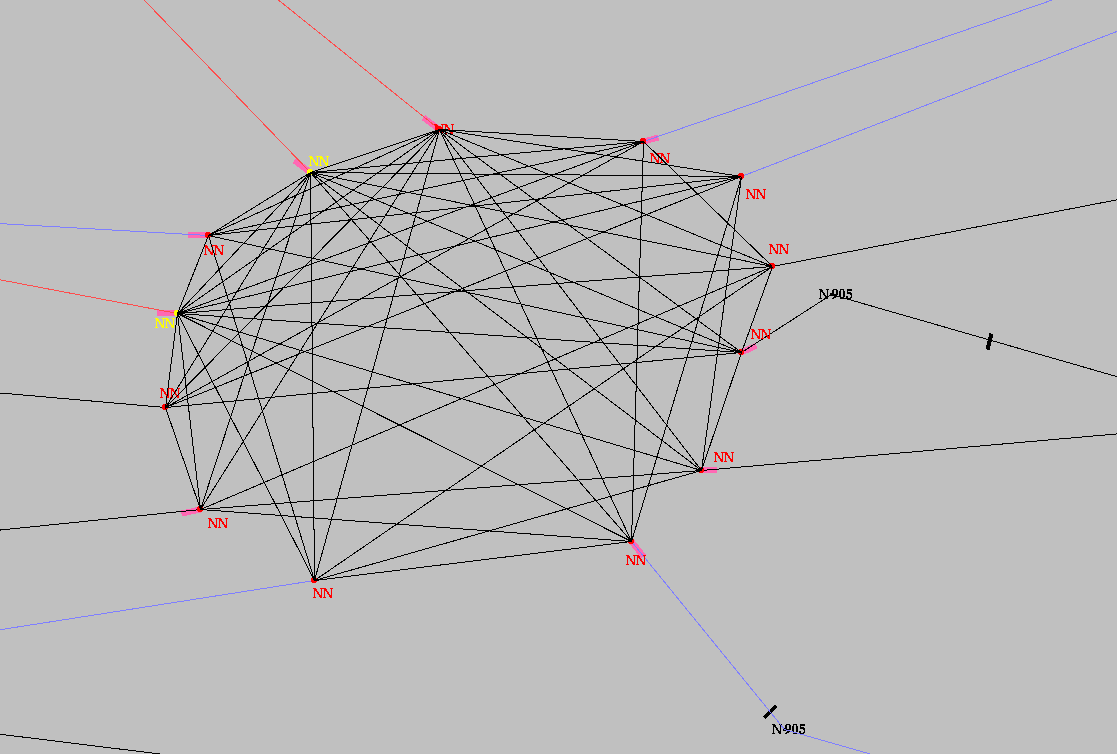}};
	\node[draw=cyan, line width=1pt, inner sep=1.5pt, fill=cyan] (NNbst) at(0.7,0.05) {};
	\node[draw=cyan, line width=1pt, inner sep=4.25pt, fill=none] (NN) at(0.575,-5.525) {};
	\draw[-latex, line width=1pt, cyan] (NN) -- (NNzoom);

	\node[draw=cyan, line width=1pt, inner sep=1.5pt, fill=cyan] at(2.25,2.75) {};
	\node[cyan] at(3.25,2.75) {Nuremberg};
	
	\node[] at(-4.25,2.75) {\footnotesize{\textbf{(a)}}};
	\node[] at(-4.25,-2.8) {\footnotesize{\textbf{(b)}}};
	\node[] at(0.75,-1.85) {\footnotesize{\textbf{(c)}}};
	
\end{tikzpicture}

%% file: table/results_complexity_routes_total.tex
\begin{tabular}{c|ccccc} 
	\toprule 
	\diagbox{$n$}{$k$} & 100 & 300 & 500 & 750 & $\infty$ \\
	\midrule 
	40 & 1148 & 2346 & 3056 & 4602 & 5916 \\
	50 & 1279 & 2479 & 3679 & 5177 & 6443 \\
 	60 & 1367 & 2567 & 3767 & 5267 & 6533 \\
	\bottomrule 
\end{tabular}

%% file: table/results_complexity_train_paths.tex
\begin{tabular}{c|ccccc} 
	\toprule 
	\diagbox{$n$}{$k$} & 100 & 300 & 500 & 750 & $\infty$ \\
	\midrule 
	40 & 63.86 & 59.18 & 60.98 & 62.2 & 58.82 \\
	50 & 98.96 & 112.30 & 115.44 & 105.8 & 118.36 \\
 	60 & 151.98 & 156.44 & 146.42 & 145.48 & 156.86 \\
	\bottomrule 
\end{tabular}

%% file: table/results_complexity_cliques.tex
\begin{tabular}{c|ccccc} 
	\toprule 
	\diagbox{$n$}{$k$} & 100 & 300 & 500 & 750 & $\infty$ \\
	\midrule 
	40 & 41.86 / 3.63 & 31.78 / 3.31 & 36.88 / 3.51 & 39.06 / 3.80 & 42.22 / 3.43 \\
	50 & 100.46 / 4.45 & 141.30 / 5.14 & 150.70 / 5.67 & 189.46 / 5.31 & 207.52 / 5.87 \\
 	60 & 316.58 / 5.83 & 250.36 / 6.67 & 236.94 / 5.99 & 480.64 / 6.39 & 475.26 / 6.34 \\
	\bottomrule 
\end{tabular}

%% file: table/results_complexity_min_max_time.tex
\begin{tabular}{c|ccccc} 
	\toprule 
	\diagbox{$n$}{$k$} & 100 & 300 & 500 & 750 & $\infty$ \\
	\midrule 
	40 & 2.6 / 0 / 14 & 2.06 / 0 / 10 & 2.32 / 0 / 14 & 2.63 / 1 / 20 & 3.86 / 1 / 105 \\
	50 & 11.8 / 1 / 88 & 15.65 / 2 / 147 & 17.54 / 0 / 158 & 16.56 / 1 / 198 & 41.08 / 4 / 1187 \\
	60 & 63.66 / 2 / 1237 & 47.38 / 1 / 376 & 47.6 / 1 / 550 & 84.16 / 0 / 1316 & 96.62 / 3 / 1210 \\
	\bottomrule 
\end{tabular}

%% file: img/time_clique_detection.tex
\begin{tikzpicture}[scale=1.0]
\begin{axis}[
    xlabel={CG iteration},
    ylabel={CPU time [s]},
    xmin=0, xmax=24,
    ymin=0, ymax=200,
    xtick={0,2,4,6,8,10,12,14,16,18,20,22,24},
    ytick={0,40,80,120,160,200},
    legend pos=north west,
    ymajorgrids=true,
    grid style=dashed,
]

\addplot[
    color=red,
    mark=otimes*,
    ]
    coordinates {
    (1,8)(2,9)(3,11)(4,15)(5,12)(6,15)(7,25)(8,18)(9,23)(10,20)(11,30)(12,29)(13,47)(14,38)(15,56)(16,51)(17,56)(18,71)(19,99)(20,97)(21,109)(22,174)(23,197)
    };
    \addlegendentry{total}  

\addplot[
    color=blue,
    mark=square*,
    ]
    coordinates {
    (1,0)(2,0)(3,0)(4,1)(5,0)(6,0)(7,1)(8,1)(9,4)(10,3)(11,4)(12,7)(13,14)(14,8)(15,23)(16,22)(17,21)(18,26)(19,41)(20,51)(21,72)(22,89)(23,115)    
    };
    \addlegendentry{clique update} 
         
\addplot[
    color=mygreen,
    mark=triangle*,
    ]
    coordinates {
    (1,8)(2,9)(3,11)(4,14)(5,12)(6,15)(7,24)(8,17)(9,19)(10,17)(11,26)(12,22)(13,33)(14,30)(15,33)(16,29)(17,35)(18,45)(19,58)(20,46)(21,37)(22,85)(23,82)
    };
    \addlegendentry{subproblem} 
\end{axis}
\end{tikzpicture}

%% file: table/results_quality.tex
\begin{tabular}{cc|cccccc} 
	\toprule 
    \multicolumn{2}{c|}{NN-$60$-$\infty$} & delay quotient & CPU time [s] & GAP [\%] & integer [\%] & \#cliques & \#iterations \\ 
    \midrule 
	gap = 0, & max CPU time = 60s & 47.53 & 18.7 & 2.0 & 92.0 & 152.9 & 6.96\\
	gap = 0.1, & max CPU time = 60s & 49.80 & 16.26 & 20.0 & 92.0 & 178.74 & 5.06 \\
	gap = 0, & max CPU time = 60s & 47.31 & 96.62 & 0.0 & 98.0 & 475.26 & 9.62 \\
	\bottomrule 
\end{tabular}

%% file: TD-GCM_revision.bbl
\begin{thebibliography}{informs2014}
\bibitem[Bettinelli et~al.(2017)Bettinelli, Santini, and Vigo]{Bettinelli2017}
A.~Bettinelli, A.~Santini, and D.~Vigo.
\newblock {A real-time conflict solution algorithm for the train rescheduling
  problem}.
\newblock \emph{Transportation Research Part B: Methodological}, 106:\penalty0
  237--265, 2017.
\newblock ISSN 01912615.
\newblock \doi{10.1016/j.trb.2017.10.005}.

\bibitem[Bonami et~al.(2015)Bonami, Lodi, Tramontani, and Wiese]{Bonami2015}
P.~Bonami, A.~Lodi, A.~Tramontani, and S.~Wiese.
\newblock {On mathematical programming with indicator constraints}.
\newblock \emph{Mathematical Programming}, 151\penalty0 (1):\penalty0 191--223,
  2015.
\newblock ISSN 14364646.
\newblock \doi{10.1007/s10107-015-0891-4}.

\bibitem[Bornd{\"{o}}rfer and Schlechte(2007)]{Borndorfer2007}
R.~Bornd{\"{o}}rfer and T.~Schlechte.
\newblock {Models for railway track allocation}.
\newblock \emph{OpenAccess Series in Informatics}, 7:\penalty0 62--78, 2007.
\newblock ISSN 21906807.

\bibitem[Bornd{\"{o}}rfer et~al.(2010)Bornd{\"{o}}rfer, Schlechte, and
  Weider]{Borndorfer2010}
R.~Bornd{\"{o}}rfer, T.~Schlechte, and S.~Weider.
\newblock {Railway track allocation by rapid branching}.
\newblock \emph{OpenAccess Series in Informatics}, 14\penalty0
  (August):\penalty0 13--23, 2010.
\newblock ISSN 21906807.
\newblock \doi{10.4230/OASIcs.ATMOS.2010.13}.

\bibitem[Br{\"{a}}nnlund et~al.(1998)Br{\"{a}}nnlund, Lindberg, N{\~{o}}u, and
  Nilsson]{Brannlund1998}
U.~Br{\"{a}}nnlund, P.~O. Lindberg, A.~N{\~{o}}u, and J.~E. Nilsson.
\newblock {Railway timetabling using Lagrangian relaxation}.
\newblock \emph{Transportation Science}, 32\penalty0 (4):\penalty0 358--369,
  1998.
\newblock ISSN 00411655.
\newblock \doi{10.1287/trsc.32.4.358}.

\bibitem[Bron and Kerbosch(1973)]{Bron1973}
C.~Bron and J.~Kerbosch.
\newblock {Algorithm 457: Finding All Cliques of an Undirected Graph}.
\newblock \emph{Communications of the ACM}, 16\penalty0 (9):\penalty0 575--577,
  1973.
\newblock ISSN 15577317.
\newblock \doi{10.1145/362342.362367}.

\bibitem[Cacchiani et~al.(2008)Cacchiani, Caprara, and Toth]{Cacchiani2008}
V.~Cacchiani, A.~Caprara, and P.~Toth.
\newblock {A column generation approach to train timetabling on a corridor}.
\newblock \emph{4or}, 6\penalty0 (2):\penalty0 125--142, 2008.
\newblock ISSN 16142411.
\newblock \doi{10.1007/s10288-007-0037-5}.

\bibitem[Cacchiani et~al.(2012)Cacchiani, Caprara, and
  Fischetti]{Cacchiani2012}
V.~Cacchiani, A.~Caprara, and M.~Fischetti.
\newblock {A lagrangian heuristic for robustness, with an application to train
  timetabling}.
\newblock \emph{Transportation Science}, 46\penalty0 (1):\penalty0 124--133,
  2012.
\newblock ISSN 15265447.
\newblock \doi{10.1287/trsc.1110.0378}.

\bibitem[Cacchiani et~al.(2014)Cacchiani, Huisman, Kidd, Kroon, Toth,
  Veelenturf, and Wagenaar]{Cacchiani2014}
V.~Cacchiani, D.~Huisman, M.~Kidd, L.~Kroon, P.~Toth, L.~Veelenturf, and
  J.~Wagenaar.
\newblock {An overview of recovery models and algorithms for real-time railway
  rescheduling}.
\newblock \emph{Transportation Research Part B: Methodological}, 63:\penalty0
  15--37, 2014.
\newblock ISSN 01912615.
\newblock \doi{10.1016/j.trb.2014.01.009}.

\bibitem[Caprara(2010)]{Caprara2010}
A.~Caprara.
\newblock {Almost 20 Years of Combinatorial Optimization for Railway Planning:
  from Lagrangian Relaxation to Column Generation}.
\newblock In \emph{10th Workshop on Algorithmic Approaches for Transportation
  Modelling, Optimization, and Systems (ATMOS'10)}, OpenAccess Series in
  Informatics (OASIcs), 2010.

\bibitem[Caprara et~al.(2002)Caprara, Fischetti, and Toth]{Caprara2002}
A.~Caprara, M.~Fischetti, and P.~Toth.
\newblock {Modeling and solving the train timetabling problem}.
\newblock \emph{Operations Research}, 50\penalty0 (5), 2002.
\newblock ISSN 0030364X.
\newblock \doi{10.1287/opre.50.5.851.362}.

\bibitem[Caprara et~al.(2006)Caprara, Monaci, Toth, and Guida]{Caprara2006}
A.~Caprara, M.~Monaci, P.~Toth, and P.~L. Guida.
\newblock {A Lagrangian heuristic algorithm for a real-world train timetabling
  problem}.
\newblock \emph{Discrete Applied Mathematics}, 154\penalty0 (5 SPEC.
  ISS.):\penalty0 738--753, 2006.
\newblock ISSN 0166218X.
\newblock \doi{10.1016/j.dam.2005.05.026}.

\bibitem[Corman et~al.(2010)Corman, D'Ariano, Pacciarelli, and
  Pranzo]{Corman2010a}
F.~Corman, A.~D'Ariano, D.~Pacciarelli, and M.~Pranzo.
\newblock {A tabu search algorithm for rerouting trains during rail
  operations}.
\newblock \emph{Transportation Research Part B: Methodological}, 44\penalty0
  (1):\penalty0 175--192, 2010.
\newblock ISSN 01912615.
\newblock \doi{10.1016/j.trb.2009.05.004}.

\bibitem[Dahms et~al.(2019)Dahms, Frank, K{\"{u}}hn, and
  P{\"{o}}hle]{Dahms2019}
F.~H.~W. Dahms, A.-L. Frank, S.~K{\"{u}}hn, and D.~P{\"{o}}hle.
\newblock {Transforming Automatic Scheduling in a Working Application for a
  Railway Infrastructure Manager}.
\newblock \emph{Link{\"{o}}ping Electronic Conference Proceedings}, 69\penalty0
  (19):\penalty0 280--289, 2019.

\bibitem[D'Ariano et~al.(2007)D'Ariano, Pacciarelli, and Pranzo]{DAriano2007a}
A.~D'Ariano, D.~Pacciarelli, and M.~Pranzo.
\newblock {A branch and bound algorithm for scheduling trains in a railway
  network}.
\newblock \emph{European Journal of Operational Research}, 183\penalty0
  (2):\penalty0 643--657, 2007.
\newblock ISSN 03772217.
\newblock \doi{10.1016/j.ejor.2006.10.034}.

\bibitem[D'Ariano et~al.(2008)D'Ariano, Corman, Pacciarelli, and
  Pranzo]{DAriano2008}
A.~D'Ariano, F.~Corman, D.~Pacciarelli, and M.~Pranzo.
\newblock {Reordering and local rerouting strategies to manage train traffic in
  real time}.
\newblock \emph{Transportation Science}, 42\penalty0 (4):\penalty0 405--419,
  2008.
\newblock ISSN 15265447.
\newblock \doi{10.1287/trsc.1080.0247}.

\bibitem[DB()]{DB}
DB.
\newblock Stammdatenaktualisierung des trassenportals (tpn).
\newblock URL
  \url{https://www.dbinfrago.com/web/aktuelles/kund-inneninformationen/kund-inneninformationen/2023-KW24-Stammdatenaktualisierung-Trassenportal-12390678}.
\newblock Accessed January 19, 2024.

\bibitem[Desaulniers et~al.(1999)Desaulniers, Desrosiers, and
  Solomon]{Desaulniers1999}
G.~Desaulniers, J.~Desrosiers, and M.~Solomon.
\newblock {Strategies for the parallel implementation of metaheuristics}.
\newblock \emph{Operations Research/ Computer Science Interfaces Series},
  15\penalty0 (January 1999):\penalty0 263--308, 1999.
\newblock ISSN 1387666X.
\newblock \doi{10.1007/978-1-4615-1507-4}.

\bibitem[Desrosiers and L{\"{u}}bbecke(2005)]{Luebbecke2005}
J.~Desrosiers and M.~E. L{\"{u}}bbecke.
\newblock {A Primer in Column Generation}.
\newblock In G.~Desaulniers, J.~Desrosiers, and M.~M. Solomon, editors,
  \emph{Column Generation}, chapter~1, pages 1--32. Springer New York, NY,
  2005.
\newblock ISBN 978-0-387-25486-9.
\newblock \doi{https://doi.org/10.1007/b135457}.

\bibitem[Dyer and Wolsey(1990)]{Dyer1990}
M.~E. Dyer and L.~A. Wolsey.
\newblock {Formulating the single machine sequencing problem with release dates
  as a mixed integer program}.
\newblock \emph{Discrete Applied Mathematics}, 26\penalty0 (2-3):\penalty0
  255--270, 1990.
\newblock ISSN 0166218X.
\newblock \doi{10.1016/0166-218X(90)90104-K}.

\bibitem[Fang et~al.(2015)Fang, Yang, and Yao]{Fang2015}
W.~Fang, S.~Yang, and X.~Yao.
\newblock {A Survey on Problem Models and Solution Approaches to Rescheduling
  in Railway Networks}.
\newblock \emph{IEEE Transactions on Intelligent Transportation Systems},
  16\penalty0 (6):\penalty0 2997--3016, 2015.
\newblock ISSN 15249050.
\newblock \doi{10.1109/TITS.2015.2446985}.

\bibitem[Lamorgese and Mannino(2012)]{Lamorgese2012}
L.~Lamorgese and C.~Mannino.
\newblock {An exact decomposition approach for the optimal real-time train
  rescheduling problem}.
\newblock \emph{PATAT 2012 - Proceedings of the 9th International Conference on
  the Practice and Theory of Automated Timetabling}, \penalty0 (October
  2023):\penalty0 428--432, 2012.

\bibitem[Lamorgese et~al.(2016)Lamorgese, Mannino, and
  Piacentini]{Lamorgese2016}
L.~Lamorgese, C.~Mannino, and M.~Piacentini.
\newblock {Optimal train dispatching by Benders'-like reformulation}.
\newblock \emph{Transportation Science}, 50\penalty0 (3):\penalty0 910--925,
  2016.
\newblock ISSN 15265447.
\newblock \doi{10.1287/trsc.2015.0605}.

\bibitem[Leutwiler and Corman(2023)]{Leutwiler2023}
F.~Leutwiler and F.~Corman.
\newblock {A review of principles and methods to decompose large-scale railway
  scheduling problems}.
\newblock \emph{EURO Journal on Transportation and Logistics}, 12\penalty0
  (December 2021), 2023.
\newblock ISSN 21924384.
\newblock \doi{10.1016/j.ejtl.2023.100107}.

\bibitem[Lodi(2010)]{Lodi2010}
A.~Lodi.
\newblock {Mixed Integer Programming Computation}.
\newblock In M.~J{\"{u}}nger, D.~Naddef, W.~R. Pulleyblank, G.~Rinaldi, T.~M.
  Liebling, G.~L. Nemhauser, G.~Reinelt, and L.~A. Wolsey, editors, \emph{50
  Years of Integer Programming 1958-2008: From the Early Years to the
  State-of-the-Art}, chapter~16, pages 619--645. 2010.
\newblock ISBN 9783540682745.
\newblock \doi{10.1007/978-3-540-68279-0}.

\bibitem[L{\"{u}}bbecke(2011)]{Luebbecke2011}
M.~E. L{\"{u}}bbecke.
\newblock {Column Generation}.
\newblock In \emph{Wiley Encyclopedia of Operations Research and Management
  Science}. 2011.
\newblock ISBN 9780470400531.
\newblock \doi{https://doi.org/10.1002/9780470400531.eorms0158}.

\bibitem[Lusby et~al.(2013)Lusby, Larsen, Ehrgott, and Ryan]{Lusby2013}
R.~M. Lusby, J.~Larsen, M.~Ehrgott, and D.~M. Ryan.
\newblock {A set packing inspired method for real-time junction train routing}.
\newblock \emph{Computers and Operations Research}, 40\penalty0 (3):\penalty0
  713--724, 2013.
\newblock ISSN 03050548.
\newblock \doi{10.1016/j.cor.2011.12.004}.

\bibitem[Maher and R{\"{o}}nnberg(2023)]{Maher2023}
S.~J. Maher and E.~R{\"{o}}nnberg.
\newblock {Integer programming column generation: accelerating branch-and-price
  using a novel pricing scheme for finding high-quality solutions in set
  covering, packing, and partitioning problems}.
\newblock \emph{Mathematical Programming Computation}, 15\penalty0
  (3):\penalty0 509--548, 2023.
\newblock ISSN 18672957.
\newblock \doi{10.1007/s12532-023-00240-w}.

\bibitem[Mascis and Pacciarelli(2002)]{Mascis2002}
A.~Mascis and D.~Pacciarelli.
\newblock {Job-shop scheduling with blocking and no-wait constraints}.
\newblock \emph{European Journal of Operational Research}, 143\penalty0
  (3):\penalty0 498--517, 2002.
\newblock ISSN 03772217.
\newblock \doi{10.1016/S0377-2217(01)00338-1}.

\bibitem[Meng and Zhou(2014)]{Meng2014}
L.~Meng and X.~Zhou.
\newblock {Simultaneous train rerouting and rescheduling on an N-track network:
  A model reformulation with network-based cumulative flow variables}.
\newblock \emph{Transportation Research Part B: Methodological}, 67:\penalty0
  208--234, 2014.
\newblock ISSN 01912615.
\newblock \doi{10.1016/j.trb.2014.05.005}.

\bibitem[Nachtigall and Opitz(2016)]{Nachtigall2016}
K.~Nachtigall and J.~Opitz.
\newblock Modelling and solving a train path assignment model.
\newblock In M.~L{\"u}bbecke, A.~Koster, P.~Letmathe, R.~Madlener, B.~Peis, and
  G.~Walther, editors, \emph{Operations Research Proceedings 2014}, pages
  423--428, Cham, 2016. Springer International Publishing.
\newblock ISBN 978-3-319-28697-6.

\bibitem[Pachl(2021)]{Pachl2021}
J.~Pachl.
\newblock \emph{{Railway Signalling Principles}}.
\newblock 2 edition, 2021.
\newblock \doi{10.24355/dbbs.084-202110181429-0}.

\bibitem[Pellegrini et~al.(2012)Pellegrini, Marli{\`{e}}re, and
  Rodriguez]{Pellegrini2012}
P.~Pellegrini, G.~Marli{\`{e}}re, and J.~Rodriguez.
\newblock {Real time railway traffic management modeling track-circuits}.
\newblock \emph{OpenAccess Series in Informatics}, 25:\penalty0 23--34, 2012.
\newblock ISSN 21906807.
\newblock \doi{10.4230/OASIcs.ATMOS.2012.23}.

\bibitem[Pellegrini et~al.(2014)Pellegrini, Marli{\`{e}}re, and
  Rodriguez]{Pellegrini2014}
P.~Pellegrini, G.~Marli{\`{e}}re, and J.~Rodriguez.
\newblock {Optimal train routing and scheduling for managing traffic
  perturbations in complex junctions}.
\newblock \emph{Transportation Research Part B: Methodological}, 59:\penalty0
  58--80, 2014.
\newblock ISSN 01912615.
\newblock \doi{10.1016/j.trb.2013.10.013}.

\bibitem[Pellegrini et~al.(2015)Pellegrini, Marli{\`{e}}re, Pesenti, and
  Rodriguez]{Pellegrini2015}
P.~Pellegrini, G.~Marli{\`{e}}re, R.~Pesenti, and J.~Rodriguez.
\newblock {RECIFE-MILP: An Effective MILP-Based Heuristic for the Real-Time
  Railway Traffic Management Problem}.
\newblock \emph{IEEE Transactions on Intelligent Transportation Systems},
  16\penalty0 (5):\penalty0 2609--2619, 2015.
\newblock ISSN 15249050.
\newblock \doi{10.1109/TITS.2015.2414294}.

\bibitem[Reynolds and Maher(2022)]{Reynolds2022}
E.~Reynolds and S.~J. Maher.
\newblock {A data-driven, variable-speed model for the train timetable
  rescheduling problem}.
\newblock \emph{Computers and Operations Research}, 142\penalty0 (February
  2021):\penalty0 105719, 2022.
\newblock ISSN 03050548.
\newblock \doi{10.1016/j.cor.2022.105719}.

\bibitem[Reynolds et~al.(2020)Reynolds, Ehrgott, Maher, Patman, and
  Wang]{Reynolds2020}
E.~Reynolds, M.~Ehrgott, S.~J. Maher, A.~Patman, and J.~Y.~T. Wang.
\newblock {A multicommodity flow model for rerouting and retiming trains in
  real-time to reduce reactionary delay in complex station areas}.
\newblock \emph{Optimization Online}, pages 1--37, 2020.

\bibitem[Sam{\`{a}} et~al.(2017)Sam{\`{a}}, D׳Ariano, Corman, and
  Pacciarelli]{Sama2017}
M.~Sam{\`{a}}, A.~D׳Ariano, F.~Corman, and D.~Pacciarelli.
\newblock {A variable neighbourhood search for fast train scheduling and
  routing during disturbed railway traffic situations}.
\newblock \emph{Computers and Operations Research}, 78:\penalty0 480--499,
  2017.
\newblock ISSN 03050548.
\newblock \doi{10.1016/j.cor.2016.02.008}.

\bibitem[Schwanh{\"{a}}u{\ss}er(1974)]{Schwanhaeusser1974}
W.~Schwanh{\"{a}}u{\ss}er.
\newblock \emph{{Die Bemessung der Pufferzeiten im Fahrplangef{\"{u}}ge der
  Eisenbahn}}.
\newblock doctoral dissertation, Fakult{\"{a}}t f{\"{u}}r Bauwesen der
  Rheinisch-Westf{\"{a}}lischen Technischen Hochschule Aachen, 1974.

\bibitem[Tomita(2017)]{Tomita2017}
E.~Tomita.
\newblock {Efficient Algorithms for Finding Maximum and Maximal Cliques and
  Their Applications}.
\newblock In S.~H. Poon, M.~S. Rahman, and H.~C. Yen, editors, \emph{WALCOM:
  Algorithms and Computation}, pages 3--15, 2017.
\newblock ISBN 9783319539249.
\newblock \doi{10.1007/978-3-319-53925-6}.

\bibitem[VIA-Con()]{Luks}
VIA-Con.
\newblock About luks.
\newblock URL \url{https://www.via-con.de/en/about-luks/}.
\newblock Accessed December 13, 2023.

\end{thebibliography}
